\documentclass{compositio}

\usepackage{isolatin1}

\usepackage{amsthm}

\usepackage{amsmath}

\usepackage{amssymb}

\usepackage{amscd}

\usepackage[all]{xy}

\usepackage{times}

\usepackage{mathrsfs}

\usepackage{mathtools}

\renewcommand{\thesection}{\arabic{section}}

\renewcommand{\theequation}{\arabic{equation}}

\theoremstyle{plain}

\newtheorem*{theorem*}{Theorem}
\newtheorem{conjecture}{Conjecture}[section]

\newtheorem{sublemma}[subsubsection]{Lemma}
\newtheorem{subtheorem}[subsubsection]{Theorem}
\newtheorem{subdefinition}[subsubsection]{Definition}
\newtheorem{subproposition}[subsubsection]{Proposition}
\newtheorem{subcorollary}[subsubsection]{Corollary}
\newtheorem{subconjecture}[subsection]{Conjecture}

\theoremstyle{definition}

\theoremstyle{remark}

\newtheorem{subremark}[subsubsection]{Remark}

\newenvironment{proofarg}[1]{\vspace{0.3cm}\noindent
\emph{Proof\nobreakspace#1.}}{\hfill
$\square$ \vspace{0.3cm}}

\numberwithin{equation}{subsubsection}

\renewcommand{\mathcal}{\mathscr}

\renewcommand{\rm}{\mathrm}

\renewcommand{\bf}{\mathbf}

\newcommand{\C}{\mathbb{C}}

\newcommand{\R}{\mathbb{R}}
\newcommand{\Q}{\mathbb{Q}}
\newcommand{\Ql}{\mathbb{Q}_{\ell}}
\newcommand{\Qp}{\mathbb{Q}_p}
\newcommand{\Z}{\mathbb{Z}}

\newcommand{\Fp}{\mathbb{F}_p}

\DeclareMathOperator{\Hom}{Hom}

\DeclareMathOperator{\Aut}{Aut}
\DeclareMathOperator{\End}{End}

\newcommand{\G}{\bf G}
\newcommand{\Gm}{\mathbb{\G}_m}

\newcommand{\der}{\mathrm{der}}

\newcommand{\Map}{\longrightarrow}

\DeclareMathOperator{\Frob}{Frob}
\DeclareMathOperator{\Gal}{Gal} 

\DeclareMathOperator{\Ind}{Ind}

\mathchardef\mhyphen="2D

\newcommand\oQ{\overline{\Q}}

\newcommand\quash[1]{}
\DeclareMathOperator{\Res}{Res}

\newcommand{\Tr}{\mathrm{Tr}} 

\usepackage{tikz}

\newcommand{\mot}{\mathrm{mot}} 
\newcommand{\class}{\mathrm{cl}} 
\newcommand{\Newton}{\mathrm{NP}} 
\newcommand{\HodgeP}{\mathrm{HTP}} 
\newcommand{\tensorInd}{\otimes\mbox{-}\Ind} 
\newcommand{\disc}{\mathrm{disc}} 
\newcommand{\cris}{\mathrm{cris}} 
\DeclareMathOperator{\rank}{rk}
\newcommand{\calM}{\mathcal{M}}

\title{Ordinary primes in Hilbert modular varieties}
\author{Junecue Suh}
\email{jusuh@ucsc.edu}
\address{Department of Mathematics \\ 1156 High St \\ University of California, Santa Cruz \\ Santa Cruz CA 95064}

\dedication{To Professor Richard Taylor with gratitude}

\begin{document}

\begin{abstract}
 A well-known conjecture, often attributed to Serre, asserts that any motive over any number field has infinitely many ordinary reductions (in the sense that the Newton polygon coincides with the Hodge polygon).
 In the case of Hilbert modular cuspforms $f$ of parallel weight $(2, \cdots, 2)$, we show how to produce more ordinary primes by using the Sato-Tate equidistribution and combining it with 
 the Galois theory of the Hecke field.  Under the assumption of stronger forms of Sato-Tate equidistribution, we get stronger (but conditional) results.
 
 In the case of higher weights, we formulate the ordinariness conjecture for submotives of the intersection cohomology of proper algebraic varieties with motivic coefficients, and verify it for the motives
 whose $\ell$-adic Galois realisations are abelian on a finite index subgroup.  We get some results for Hilbert cuspforms of weight $(3,\cdots, 3)$, weaker than those for $(2, \cdots, 2)$.
\end{abstract}

\maketitle

\section{Introduction}
\label{sec-intro}

Let $\mathcal{X}$ be a projective smooth scheme over the ring $\mathcal{O}_{F,S}$ of $S$-integers in a number field $F$, and let $X/F$ be the generic fibre. 
For each integer $i$, one defines the Hodge polygon $HP^i(X)$ from the Hodge filtration on $H^i_{dR}(X/F)$: it is the convex planar graph emanating from the origin in which the line segment of slope $a$ appears with  the multiplicity equal to the Hodge number $h^{a, i-a} = \dim_F H^{i-a}(X, \Omega^a)$.
For every $\mathfrak{p}$ outside the finite set $S$, one also forms the Newton polygon $\Newton(\Frob_{\mathfrak{p}}, H^i_{\cris}((\mathcal{X}\otimes_{\mathcal{O}_F} k(\mathfrak{p}))/W(k(\mathfrak{p}))))$ using the $p$-adic slopes of the crystalline Frobenius.  (See \S \ref{sec-semiring-newton}.)

Katz \cite[Conjecture 2.9]{katz-ax} conjectured, and Mazur proved \cite{mazur}, that the latter lies above the former:
$$
\Newton(\Frob_{\mathfrak{p}}, H^i_{\cris}((\mathcal{X}\otimes_{\mathcal{O}_F} k(\mathfrak{p}))/W(k(\mathfrak{p}))) \ge HP^i(X).
$$
When the two coincide, we say that $\mathfrak{p}$ is an \textit{ordinary} prime (in degree $i$ for $X$), following Mazur.  
The following appears to have been first considered by Serre \cite[no. 133]{serre-oeuvres4} for abelian varieties $X$ (in degree $1$):
\begin{conjecture}[(Ordinariness Conjecture)]
\label{conj-serre}
There exists an infinite set of primes $\mathfrak{p}$ such that
$$
\Newton(\Frob_{\mathfrak{p}}, H^i_{\cris}((\mathcal{X}\otimes_{\mathcal{O}_F} k(\mathfrak{p}))/W(k(\mathfrak{p}))) = HP^i(X).
$$
\end{conjecture}
It is known to be true for elliptic curves and abelian surfaces by arguments of Katz and of Ogus \cite[Prop. 2.7]{ogus-hodge}, and for abelian varieties whose endomorphism ring is $\Z$ and whose algebraic monodromy group satisfies a condition: See Pink \cite[\S 7]{pink-mtc}.
It is also known for all CM abelian varieties. [$\ast\ast$]

In this article, we investigate Conjecture \ref{conj-serre} for the factors of modular Jacobians cut out by cusp forms of weight $2$,
and provide several methods for finding ordinary primes in them.

More generally, we consider the parts of the intersection cohomology of the Hilbert modular varieties attached to totally real number fields $F$ of degree $d=[F:\Q]$,
cut out by new normalised cusp forms $f$ of parallel weight $(2, \cdots, 2)$.  The ordinariness in this context first appeared as an assumption in the construction of Galois representations by Wiles \cite{wiles},
which was later removed by Taylor \cite{taylor-hilbert} and by Blasius and Rogawski \cite{blasius-rogawski}.  
However, many of the best results and constructions in Iwasawa and Hida theory at present depend on the existence of ordinary primes in a crucial way: See among others the works of Emerton, Pollack and Weston \cite{emerton-pollack-weston}; of Ochiai \cite{ochiai}; of Nekov\'{a}\v{r} \cite{nekovar-parity-2}; of Dimitrov \cite{dimitrov-ajm}; of Skinner and Urban \cite{skinner-urban}; and of Wan \cite{wan}.

Since, at the moment, we lack a satisfactory `crystalline' theory of perverse sheaves or intersection cohomology,
we first formulate in \S \ref{sec-form-ih} analogues of the Katz Conjecture and Conjecture \ref{conj-serre} for the $\ell$-adic \'etale intersection cohomology of any projective variety $X$,
where $\ell$ is any auxiliary prime.
Here we form the \textit{Hodge-Tate polygon} attached to the Galois representation in place of the Hodge polygon, and the Newton polygon by using the $\ell$-adic Frobenius.

These conjectures (Conjectures \ref{conj-ih-katz} and \ref{conj-ih-serre}) satisfy a basic consistency, in that (a) the Hodge-Tate polygon is independent of the choice of a prime $\lambda$ of $F$ lying over $\ell$; 
and (b) the conjectures are independent of the auxiliary prime $\ell$. 
We prove these statements by using theorems of Gabber (on the independence of $\ell$ in the intersection cohomology of complete varieties) \cite{gabber}, of Katz and Laumon (on the constructibility properties of certain constructions in derived categories) \cite{katz-laumon}, of Andr\'e (his theory of motivated cycles) \cite{andre-motifs} and of de Cataldo and Migliorini (on the motivated nature of the decomposition theorem in intersection cohomology) \cite{decataldo-migliorini}.

Let us return to the Hilbert modular varieties and the forms $f$.
By using either (i) the theorem of de Cataldo and Migliorini in \textit{op. cit.} and its rational extension due to Patrikis \cite{patrikis} or (ii) the recent motivic constructions of Ivorra and Morel \cite{ivorra-morel}, we construct an intersection motive of $X$ which in realisations give the intersection cohomology. 
Then by lifting the action of the Hecke correspondences on the intersection cohomology to one on the intersection motive (see Proposition \ref{rmk-hecke-andre-motives}), we construct an Andr\'e motive $M(f)$.
The conjectures make sense for these submotives, and the consistency mentioned above also hold for them.  We denote by $K_f$ the Hecke field of $f$.


We say that a subset $\Sigma$ of $\mathrm{MaxSpec}(\mathcal{O}_F)$ is \textit{abundant} if $\Sigma$ has lower (natural) density $>0$,
and that $\Sigma$ is \textit{principally abundant} if there exists a finite extension $F'/F$ such that the inverse image of $\Sigma$ in $\mathrm{MaxSpec}(\mathcal{O}_{F'})$
has density $=1$ in $F'$.
In the previous cases where Conjecture \ref{conj-serre} has been established, in fact a principally abundant set of ordinary primes was found.

By using the construction of the Galois representation attached to $f$ and the purity of $IH$, we first show that $M(f)$ satisfies the analogue of the Katz Conjecture
and that we can push the Newton polygon `half way' to ordinariness in a quantifiable sense, for a principally abundant set of primes.
However, in the attempt to push just beyond the half-way threshold, we face an obstruction of ``geometry of numbers'' type (Minkowski).
We show how to overcome it (for a principally abundant set of primes) by using a stronger form of Sato-Tate equidistribution (see \S \ref{sec-sst}), but this last form remains unknown in general.

In order to go further and to obtain unconditional results, we look into (a) `multivariate' variants of the Sato-Tate Conjecture in \S \ref{sec-sst} and (b) the interaction between $F$ and $K_f$ in \S \ref{sec-interplay}. 
For the latter, we define an invariant, the \textit{slope} $\sigma_F(K) \in [0, 1]$ of a coefficient number field $K$ over a ground number field $F$ (see Definition \ref{def-slope}),
by using the action of $\Gal(\oQ/F)$ on the set $\Hom(K, \oQ)$ of field embeddings.
The slope $\sigma_F(K)$ is $0$ if, for example, (1) $[K:\Q]$ is a prime not dividing $[F:\Q]$; or 
(2) the Galois group of $K/\Q$ is the (full) symmetric group on $[K:\Q]$ letters\footnote{this `Maeda-like' condition appears to be often satisfied for the Hecke fields $K = K_f$ in practice, but not always.  See Section 5 for examples.} and $K$ and $F$ have coprime discriminants.

Here is a collection of results which follow from the Main Theorem \ref{thm-main}:
\begin{theorem*}
Let the notation be as above. Then $M(f)$ has an abundant set of ordinary primes if at least one of the following conditions is satisfied:
\begin{enumerate}
 \item[(a)] $[K_f^{\circ}:\Q]\le 2$; where $K_f^{\circ}$ is the smallest Frobenius field of $f$ in the sense of Ribet (see \S \ref{sec-ribet} for the precise definition);
 \item[(b)] $f$ is of CM type;
 \item[(c)] the slope $\sigma_{\widetilde{F}}(K_f)$ is equal to $0$, where $\widetilde{F}$ is the Galois closure of $F$ over $\Q$; or
 \item[(d)] an element of $\Gal(\oQ / \widetilde{F})$ has exactly $2$ orbits in $\Hom(K_f, \oQ)$, the orbits have the same size 
 and $f$ satisfies a strong form of Sato-Tate equidistribution named (RST) in \S \ref{sec-sst}.
\end{enumerate}
\end{theorem*}

The statement in (c) follows from the quantitative part (4) of Theorem \ref{thm-main}: The smaller the slope $\sigma$, the closer to ordinariness we can push the Newton polygons.

In Section \ref{sec-numerical}, we look at the forms $f$ with low levels for $4$ number fields $F$ of degree $d\le 4$,
and show that for most $f$ under consideration, (a), (b) and (c) provide an abundant set of ordinary primes unconditionally, and that (d) complements them under the strong Sato-Tate condition.
We give descriptions of the conditional cases as well as some cases where our methods fall short of yielding abundance of ordinary primes.

\noindent\textbf{Plan.} Here are some main ingredients and ideas in the text.

In Section 2, we show (Theorem \ref{thm-padic-complex}) that for Andr\'e motives, the Hodge polygon made of the Hodge numbers in its (transcendental) Betti realisations coincide with the Hodge-Tate polygon made of the $p$-adic Hodge-Tate weights in the (algebraic) $p$-adic \'etale realisations.
This applies in particular to the Hecke isotypic components of the intersection cohomology (motive) of the Hilbert modular varieties.

In Section 3, we introduce a few things in preparation for the Main Theorem in Section 4. (a) In \S \ref{sec-semiring-newton} we study the partially ordered semiring of Newton polygons; this will be useful in dealing with the tensor induction in Theorem \ref{thm-main}. (b) In \S \ref{sec-interplay}, we define the notion of \emph{slope} and \emph{bisection} in the interaction between the ground field $F$ and the coefficient field $K$; very roughly speaking, they measure the sizes and shapes of `large' orbits in the Hecke and Frobenius fields.  In Theorem \ref{thm-main}, these will be combined with the Chebotarev density theorem and strong forms of the Sato-Tate Conjecture.

In Section 4 we prove the Main Theorem \ref{thm-main}. Here we make a connection (which to the author's knowledge is new) between two conjectures on the eigenvalues of the Frobenius elements: Namely, the Sato-Tate Conjecture on the archimedean properties of the eigenvalues (which in turn is closely linked to the Langlands Program) on the one hand and the Ordinariness Conjecture on the $p$-adic properties (with varying $p$) of the eigenvalues on the other hand.

In the final Section \ref{sec-gen-mot-coeffs}, we formulate analogues of the Katz and the ordinariness conjectures for submotives of the intersection cohomology of more general motivic coefficients, following a suggestion of Katz.
In cases where we have good crystalline realisations compatible with the $\ell$-adic realisations (which include the nonconstant motivic coefficients on Hilbert modular varieties), 
we verify the Katz conjecture by using Mazur's theorem.  In case the submotive has potentially abelian $\ell$-adic realisation, we also verify the ordinariness conjecture by using Serre's theory \cite{serre-mcgill}.
Finally, we provide some methods to deal with the parallel motivic weight $(3, \cdots, 3)$ in the Hilbert modular case.

\section{Formulation of Conjectures for $IH$}
\label{sec-form-ih}

\subsection{Polygons}

Let $F$ be a number field with algebraic closure $F^s$, $\ell$ a prime number, $V$ a $\Ql$-vectorspace of dimension $m<\infty$, and
$$
\rho: \Gal(F^s / F) \Map \Aut(V) \simeq GL_m(\Ql)
$$
a continuous representation that is unramified outside a finite set $S$ of maximal ideals of $\mathcal{O}_F$.

\begin{subdefinition}
Assume that $\rho$ is $\Q$-rational in the sense of Serre \cite{serre-mcgill}.

Then for each maximal ideal $\mathfrak{p}$ of $\mathcal{O}_F$ outside $S$ and with residue characteristic $p\neq \ell$,
we define the \textbf{Newton polygon} $\Newton(\Frob_{\mathfrak{p}}, \rho)=\Newton(\Frob_{\mathfrak{p}}|_V)$
as the Newton polygon of the characteristic polynomial
$$
\det( T - \Frob_{\mathfrak{p}}: V) \in \Q[T]
$$
with respect to the $p$-adic valuation $v_{\mathfrak{p}}$ on $\Q$ normalised by $v_{\mathfrak{p}}( \mathbb{N}\mathfrak{p}) = 1$.
\end{subdefinition}
Equivalently: Choose an isomorphism $\oQ_{\ell} \simeq \oQ_p$ and let $x_1, \cdots, x_m \in \oQ_{\ell}$ be the eigenvalues of $\Frob_{\mathfrak{p}}$ on $V$.
Then, the multiset of slopes
$$
\left\{ \frac{v_p(x_1)}{v_p(\mathbb{N}\mathfrak{p})}, \cdots, \frac{v_p(x_m)}{v_p(\mathbb{N}\mathfrak{p})} \right\}
$$
gives the Newton polygon.  It is independent of the chosen isomorphism.

\begin{subdefinition}
Assume that $\rho$ is Hodge-Tate at every prime $\lambda$ of $F$ lying over $(\ell)$, and that the set with multiplicities of the Hodge-Tate weights at $\lambda$ is independent of $\lambda|\ell$.

Then we define the \textbf{Hodge-Tate polygon} $\HodgeP(\rho)=\HodgeP(V)$ as the convex planar polygon starting from $(0,0)$ in which the slope $i$ appears as many times as the Hodge-Tate weight $i$ appears in $\rho$.
\end{subdefinition}
There seem to be competing sign conventions for the Hodge-Tate weights.  We take the `geometric' one, so that $H^1$ of an elliptic curve has Hodge-Tate weights $\{0, 1\}$.

\subsection{Independence of Hodge-Tate weights}

First, we show that the Hodge-Tate polygon we have defined coincides with the classical Hodge polygon for all Andr\'e motives, thereby extending a theorem of Faltings \cite{faltings}, \cite{faltings-jami}.
For this we will crucially rely on the theory of motivated cycles and the resulting category of Andr\'e motives, given in \cite{andre-motifs}.

\begin{subtheorem}
\label{thm-padic-complex}
Let $M$ be an Andr\'e motive over a finite extension field $K$ of $\Qp$, and let $\sigma:K\Map \C$ be a complex embedding.
Denote by $M_p$ its $p$-adic \'etale realisation, and by $M_{\sigma}$ its Betti realisation via $\sigma$.

Then the set with multiplicities of the Hodge-Tate weights of $M_p$ coincides with that of the complex Hodge numbers of $M_{\sigma}$.
\end{subtheorem}

\begin{proof}
We may and will assume that $M$ is simple, and that there exist a projective smooth variety $Y$ of dimension $d$ over $K$, an integer $n$, and an Andr\'e motivated cycle
$$
\xi \in A^d_{\mot}(Y \times_K Y)
$$
such that $\xi$ acts as the idempotent cutting out $M$ in $\mathfrak{h}^n(Y)$.  
Let $e_p$ (resp. $e_{\sigma}$, resp. $e_{dR}$) be the image of $\xi$ in the $p$-adic \'etale (resp. $\sigma$-Betti, resp. de Rham) realisation:
$$
e_p \in H^{2d}((Y\times Y)\otimes_K \overline{K}, \Qp)(d), \,\, e_{\sigma} \in H^{2d}(\sigma(Y\times Y), \Q(d)), \,\, e_{dR} \in H^{2d}_{dR}(Y\times Y)(d).
$$

We have the following diagram 
$$
\xymatrix{
 & A^d_{\mot} (Y \times Y) \ar[ld]_-{\class_{dR}} \ar[d]_-{\class_{\sigma}} \ar[rd]^-{\class_p} & \\
 H_{dR} \ar[d] \ar@/_3pc/[dd] & H_{\sigma} \ar[d] \ar[dr] & H_p \ar[d]^{\mathrm{Art}}_-{\sim} \ar@/^3pc/[dd] \\
 H_{dR}\otimes_K \C & H_{\sigma}\otimes \C  \ar[l]^-{\mathrm{comp}_{\C}}_-{\sim} & H_{\sigma}\otimes \Qp \\
 H_{dR} \otimes_K B_{dR} & & H_p \otimes_{\Qp} B_{dR} \ar[ll]_-{\sim}^{\mathrm{comp}_{dR}} 
}
$$
Here we suppressed $Y\times Y$ in the argument for all cohomology theories as well as the degree $2d$ and the Tate twist $(d)$.
Undecorated arrows are extensions of scalars and Art denotes Artin's comparison isomorphism.

\noindent \textit{Main Point}: The diagram is commutative.  
In particular, the image of $\xi$ in any group in the diagram is the same, no matter which path emanating from $A^d_{\mot}(Y\times Y)$ is followed.

This follows from the definition of Andr\'e motivated cycles,
together with the fact that the comparison isomorphisms in display are isomorphisms between Weil cohomology theories (and as such
compatible with pullback, pushforward, cup product, cycle class, and Poincar\'e duality, that are involved in the definition of Andr\'e motivated cycles).
See Andr\'e \cite[\S 2.3 and \S 2.4]{andre-motifs}. 

(In contrast, it is not clear whether the similar diagram would be commutative, if we replace the apex $A_{\mot}$ with the larger space of the absolute Hodge cycles.)

Now applying the idempotents obtained from $\xi$ to the similar diagram without apex:
$$
\xymatrix{
 H^n_{dR}(Y) \ar[d] \ar@/_4pc/[dd] & H^n(\sigma(Y), \Q) \ar[d] \ar[dr] & H^n(Y\otimes_K \overline{K}, \Qp) \ar[d]^{\mathrm{Art}}_-{\sim} \ar@/^4pc/[dd] \\
 H^n_{dR}(Y) \otimes_K \C & H^n(\sigma(Y), \Q) \otimes \C  \ar[l]^-{\mathrm{comp}_{\C}}_-{\sim} & H^n(\sigma(Y), \Q) \otimes \Qp \\
 H^n_{dR}(Y) \otimes_K B_{dR} & & H^n(X\otimes_{K} \overline{K}, \Qp) \otimes_{\Qp} B_{dR} \ar[ll]_-{\sim}^{\mathrm{comp}_{dR}} 
}
$$
we get the diagram
$$
\xymatrix{
 M_{dR} \ar[d] \ar@/_3pc/[dd] & M_{\sigma} \ar[d] \ar[dr] & M_p \ar[d]^{\mathrm{Art}}_-{\sim} \ar@/^3pc/[dd] \\
 M_{dR} \otimes_K \C & M_{\sigma} \otimes \C  \ar[l]^-{\mathrm{comp}_{\C}}_-{\sim} & M_{\sigma} \otimes \Qp \\
 M_{dR} \otimes_K B_{dR} & & M_p \otimes_{\Qp} B_{dR} \ar[ll]_-{\sim}^{\mathrm{comp}_{dR}} 
}
$$

Now, on the one hand, the Hodge-Tate weights of $M_p$ can be read off from the filtration on 
$$
(M_p \otimes_{\Qp} B_{dR})^{\Gal(\overline{K} / K)} \simeq M_{dR}  \,\,\,\,\, \mbox{ (filtered isomorphism) }.
$$
On the other hand, the Hodge numbers of $M_{\sigma}$ can also be read off from the (algebraic) Hodge filtration on $M_{dR}$ in a similar manner.

\end{proof}

\begin{subcorollary}
Let $X$ be any projective variety over any finite extension $K$ of $\Qp$, and let $\sigma: K\Map \C$ be any complex embedding.
Then for any integer $n$, the Hodge-Tate weights of $IH^n(X\otimes_K \overline{K}, \Qp)$ coincide with the Hodge numbers of $IH^n(\sigma(X), \Q)$.
\end{subcorollary}

\begin{proof}
Let $\pi: Y \Map X$ be any resolution of singularities.  We use the main theorem of de Cataldo and Migliorini \cite{decataldo-migliorini}, 
strengthened in $K$-rationality by Patrikis \cite[\S 8]{patrikis} 
to deduce the existence of an Andr\'e motivated cycle
$$
\xi := \xi_n \in A^d_{\mot} (Y \times Y)
$$
such that the Betti realisation
$$
e_{\sigma}:= \class_{\sigma} (\xi) \in H^{2d}( \sigma(Y\times Y), \Q(d))
$$
defines as a correspondence the idempotent for the direct summand
$$
IH^n (\sigma(X), \Q) \subseteq H^n(\sigma(Y), \Q)
$$
and similarly for the $p$-adic \'etale realisation
$$
e_{p} := \class_{\Qp} ( \xi ) \in H^{2d} ( (Y\times Y)\otimes_K \overline{K}, \Qp (d) ).
$$

Using this, de Cataldo-Migliorini and Patrikis \textit{define} the $K$-rational intersection de Rham cohomology of $X$:
$$
IH^n_{dR}(X/K)
$$
as the image of the idempotent
$$
e_{dR} := \class_{dR}(\xi) \in H^{2d}(Y\times Y)(d)
$$
acting on $H^n_{dR}(Y)$.  Clearly, there is a comparison isomorphism of $IH^n_{dR}(X/K)\otimes_K \C$ with the (transcendental) Hodge structure on the Betti realisation $IH^n(\sigma(X), \Q)\otimes \C$.
Moreover, this last Hodge structure coincides with the Hodge structure contsructed by Morihiko Saito, see de Cataldo \cite[Th. 4.3.5]{decataldo}.

Apply Theorem \ref{thm-padic-complex} to this situation.
\end{proof}

\begin{subcorollary}
\label{cor-indep-andre}
Let $M$ be an Andr\'e motive over a number field $K$, $\mathfrak{P}$ any maximal ideal of $\mathcal{O}_K$ of residue characteristic $p$, and $\sigma: K \Map \C$ any complex embedding.

Then the Hodge-Tate weights of $M_p$ at $\mathfrak{P}$ coincide with the Hodge numbers of $M_{\sigma}$.  In particular, the Hodge-Tate weights are independent of $\mathfrak{P} \in \mathrm{MaxSpec} (\mathcal{O}_K)$.
\end{subcorollary}

This applies, for example, to the Andr\'e motives cut out by \textit{algebraic} cycles from the cohomology of a projective smooth variety over $K$.


\begin{subcorollary}
Let $X$ be a projective variety defined over a number field $K$, $\mathfrak{P}$ any maximal ideal of $\mathcal{O}_K$ of residue characteristic $p$, and $\sigma: K \Map \C$ any complex embedding.

Then for every integer $n$, the Hodge-Tate weights of $IH^n(X\otimes_K \overline{K}, \Qp)$ at $\mathfrak{P}$ coincide with the Hodge numbers of $IH^n (\sigma(X), \Q)$.

In particular, the Hodge-Tate weights are independent of $\mathfrak{P} \in \mathrm{MaxSpec} (\mathcal{O}_K)$.
\label{cor-ind-htw-ih}
\end{subcorollary}

\begin{subremark}
Strictly speaking, the results on the intersection cohomology can be proven by using the de Cataldo-Migliorini theorem \cite{decataldo-migliorini} only
(and not using the $K$-rational version \cite{patrikis}).
To see this, note that the Hodge numbers and the Hodge-Tate weights are insensitive to the base change to a finite extension of $K$ (in both local and global cases),
and the construction of \cite{decataldo-migliorini} yields the necessary Andr\'e motivated cycle over a finite extension of $K$.
\end{subremark}

In case $X=X^{BB}$ is the Baily-Borel compactification of a Shimura variety $X^{\circ}$ (we refer to Ash-Mumford-Rapoport-Tang \cite{AMRT} and Pink \cite{pink-mtc} in general, and Brylinski-Labesse \cite{BL} and Rapoport \cite{rapoport} in the special case of Hilbert modular varieties), one further decomposes the intersection cohomology of $X$ into the Hecke-isotypic components: the Hecke correspondences act on the intersection cohomology of $X$, and span a  $\Q$-subalgebra $\mathcal{H}_{X, \Q}$ in the finite dimensional $\Q$-algebra $\End_{HS}(IH^d_B (X, \Q))$.
By decomposing $\mathcal{H}_{X, \Q}$ into a product of $\Q$-algebras, we obtain the Hecke isotypic components.

\begin{subproposition}
 \label{rmk-hecke-andre-motives}
 The Hecke isotypic components come from Andr\'e motives over the reflex field $E$.  
 
 As a consequence, Theorem \ref{thm-padic-complex} applies to these components.
\end{subproposition}
\begin{proof}
 This boils down to first finding an Andr\'e (pure Nori) motive 
 $\mathfrak{ih}(X)$ whose $\ell$-adic and Betti realisations give the $\ell$-adic and Betti intersection cohomology of $X$; and then lifting the action of the Hecke correspondences to one on $\mathfrak{ih}(X)$.
 
 While the first step can be done as above in an ad hoc fashion --- using a (noncanonical) resolution of singularities and using the theorems of de Cataldo, Migliorini and Patrikis --- the second step is done most systematically (in our opinion) by using the theory of weight filtration.  We learned the argument from S. Morel (cf. \cite[\S 5]{morel-jams}), which uses the more recent motivic constructions of Ivorra and Morel \cite{ivorra-morel}.

 (Before the details, let us stress the main point and indicate where the innovations are. The intersection complex $IC$ (as a perverse sheaf in the derived category of constructible sheaves) was originally constructed as an iterated application of the $2$-step procedure: taking the direct image $R j_{\ast}$ under open immersions $j$, and then truncating with respect to a topological stratification and a function called `perversity'. See the explicit formula \cite[Prop. 2.1.11]{bbd}. 
 
 One of the innovations in \cite{morel-jams} was to realise $IC$ (for the middle perversity) over finite fields as the truncation with respect to the \emph{weight} filtration: See \cite[Th. 3.1.4]{morel-jams}.
 
 This renders the extension of the Hecke operators to $IC$ deceptively easy: One needs neither to worry about singularities in the boundary (which can be bad) --- to which one must pay heed if one uses the original (topological) definition --- nor to rely on the toroidal compactifications --- of which there is no canonical choice, and many are needed to extend the Hecke correspondences.
 This is why we adopt the idea.
 
 In \cite{ivorra-morel}, Ivorra and Morel construct the four operations of Grothendieck (namely $f^{\ast}_{\calM}$, $f^{!}_{\calM}$, $f_{\ast}^{\calM}$, $f_{!}^{\calM}$ for morphisms $f$ between quasiprojective varieties) and the weight filtrations on the derived category of `perverse mixed motives'; this last abelian category is moreover shown to be equivalent to the Nori motives.
 
 This allows one to formally apply the algorithm described in \cite[\S 5]{morel-jams}, but this time applied in the motivic derived category of \cite{ivorra-morel}, rather than in the derived category of constructible sheaves as in \cite{morel-jams}. Then the Betti and $\ell$-adic realisations of these constructions give rise to the Hecke operators constructed previously on the level of complexes of sheaves.)

 Now let us turn to a more detailed argument, and indicate which constructions in \cite{ivorra-morel} replace those in the parallel argument from \cite[\S 5]{morel-jams}.
 
 Let $j:X^{\circ} \hookrightarrow X=X^{BB}$ be the open immersion of the Shimura variety and $c_1, c_2: Y^{\circ} \Map X^{\circ}$ be the two finite \'etale maps that give rise to a given Hecke correspondence. With $j':Y^{\circ} \hookrightarrow Y$ denoting the open immersion into the Baily-Borel compactification, the $c_i$ extend canonically to $\overline{c_i}:Y\Map X$.  We start with the identity correspondence
 $$
 u=1: c_1^{\ast} \Q_{X^{\circ}} \Map c_2^{!} \Q_{X^{\circ}}.
 $$
 arising from the natural isomorphisms $c_1^{\ast}\Q_{X^{\circ}} \simeq \Q_{Y^{\circ}}$ and $c_2^{!}\Q_{X^{\circ}} \simeq \Q_{Y^{\circ}}$.
 
 Here and below, the functors $f^{\ast}$ and $f^{!}$, etc., refer to the ones  in \cite[Th. 5.1]{ivorra-morel}, where the notations $f^{\ast}_{\calM}$, $f^{!}_{\calM}$, etc., are used.
  
 Just as in \cite[\S 5.1]{morel-jams}, but using the motivic constructions of the $4$ functors and the base change morphisms, stated \cite[Th. 5.1]{ivorra-morel}, we take $Rj'_{\ast} = {j'}_{\ast}^{\calM}$ and use the base change morphisms
 $$
 \xymatrix{
 \overline{c_1}^{\ast} Rj_{\ast} \Q_{X^{\circ}}
 \Map Rj'_{\ast} c_1^{\ast} \Q_{X^{\circ}} \ar[r]^-{u} & Rj'_{\ast} c_2^{!} \Q_{X^{\circ}} \Map \overline{c_2}^{!} Rj_{\ast} \Q_{X^{\circ}}.
 }
 $$
 We then use the fact that the lowest weight filtration of $Rj_{\ast} \Q_{X^{\circ}}$ is canonically isomorphic to the intersection complex $j_{!\ast} (\Q_{X^{\circ}}[d])[-d]$.
 For this, it is enough to show that the (motivic) weight filtration has $\ell$-adic realisation equal to the ($\ell$-adic) weight filtration. But this follows from the definitions and construction: \cite[Def. 6.12, Def. 6.13, and Prop. 6.16]{ivorra-morel}.
 
 Then we proceed as in \cite[Lem. 5.1.4]{morel-jams} to obtain
 $$
 \overline{u}: \overline{c_1}^{\ast} IC_X \Map \overline{c_2}^{!} IC_X,
 $$
 lifting the cohomological correspondence in realisations.  Again the point is that all the arrows in display in \cite[Lem. 5.1.4]{morel-jams} are constructed using (only) the functoriality of $Rj_{\ast}$, the base change morphisms, and the weight filtration.  In our (motivic) context, we use the main theorem \cite[Th. 5.1]{ivorra-morel} for the first two and the weight filtration constructed in \cite[Prop. 6.16]{ivorra-morel} for the last.
\end{proof}

This applies in particular to the Hilbert modular varieties and the submotive $M(f)$ of $\mathfrak{ih}^d(X^{BB})$ cut out by any new cuspform $f$ of parallel weight $(2, \cdots, 2)$ and its conjugates,
and Corollary \ref{cor-indep-andre} applies to it.

\subsection{Conjectures}

Now we can formulate the analogue of Katz's Conjecture:
\begin{subconjecture}
Let $X$ be a projective variety over a number field $F$, and let $n$ be an integer.

Then there exists a finite set $S=S(X, n)$ of maximal ideals of $\mathcal{O}_F$ such that for every prime number $\ell$ and every maximal ideal $\mathfrak{p}$ of $\mathcal{O}_F$ outside $S$ and with residue characteristic $\neq \ell$,
we have
$$
\Newton(\Frob_{\mathfrak{p}}|_{IH^n(X\otimes_F F^s, \Ql)}) \ge \HodgeP(IH^n(X\otimes_F F^s, \Ql)).
$$
\label{conj-ih-katz}
\end{subconjecture}
In case $X$ is also smooth, this is known to be true, by theorems of Katz and Messing \cite{katz-messing} (comparing the $\ell$-adic Frobenius at $\mathfrak{p}$ with the crystalline Frobenius), of Mazur \cite{mazur} (showing that the Newton polygon lies on or above the Hodge polygon for the crystalline cohomology), and of Faltings \cite{faltings} (showing that the Hodge-Tate polygon is the same as the Hodge polygon).

And we formulate the analogue of the `ordinariness' conjecture:
\begin{subconjecture}
Let $X$ be a projective variety over a number field $F$, and let $n$ be an integer.

For every prime number $\ell$, there exists an infinite set of maximal ideals $\mathfrak{p}$ of $\mathcal{O}_F$ with residue characteristic $\neq \ell$ such that
$$
\Newton(\Frob_{\mathfrak{p}}|_{IH^n(X\otimes_F F^s, \Ql)}) = \HodgeP(IH^n(X\otimes_F F^s, \Ql)).
$$
\label{conj-ih-serre}
\end{subconjecture}

We note that the right hand side of the conjectures is independent of $\lambda$ or $\ell$ by Corollary \ref{cor-ind-htw-ih},
and that the left hand side is independent of $\ell$ because the $IH^n(X, \Ql)$ form a strictly compatible system by Gabber \cite{gabber} and Katz-Laumon \cite[Th. 3.1.2]{katz-laumon}.
Therefore Conjectures \ref{conj-ih-katz} and \ref{conj-ih-serre} are independent of the auxiliary prime $\ell$.

In case $X$ is also smooth, Conjecture \ref{conj-ih-serre} is equivalent to Conjecture \ref{conj-serre} recalled in \S \ref{sec-intro}, 
since (1) the Newton polygons of the $\ell$-adic and crystalline Frobenius endomorphisms are the same by Katz and Messing \cite{katz-messing},
and (2) the Hodge-Tate polygon coincides with the Hodge polygon by Faltings \cite{faltings}.


\section{Preparation}

\subsection{Notation}
\label{sec-notation}

From this point on, $F\subseteq \oQ$ denotes a totally real number field of degree $d=[F:\Q]$ and discriminant $\disc(F)$; $\widetilde{F}$ is the Galois closure of $F/\Q$, of degree $\widetilde{d} = [\widetilde{F}: \Q]$.

Let $f$ be a new normalised Hilbert eigencuspform of parallel weight $(2, \cdots, 2)$ of level $\mathfrak{n}\subseteq \mathcal{O}_F$. 
The Fourier coefficients of $f$ generate the number field:
$$
K_f := \Q \left( \left\{ a_{\mathfrak{p}} \right\}_{\mathfrak{p}}  \right),
$$
where $\mathfrak{p}$ ranges over the primes of $\mathcal{O}_F$ not dividing $\mathfrak{n}$.
It is either a totally real number field or a CM field, and we let $k_f:=[K_f:\Q]$.

We note that the ordinariness in this context is equivalent to the following simple condition: $\mathfrak{p}$ is an ordinary prime (for $f$) if and only if $a_{\mathfrak{p}}$ is nonzero and does not belong to any prime ideal $\wp$ of $\mathcal{O}_{K_f}$ lying over $(p) = \mathfrak{p}\cap \Z$.

We fix once and for all a rational prime $\ell$ that splits completely\footnote{we choose a split prime just for simplifying the exposition a little bit.
The obvious analogues of Conjectures \ref{conj-ih-katz} and \ref{conj-ih-serre} are independent of $\ell$ for $M(f)$ defined below, since we still have a strictly compatible system of Galois representations.} in $K_f$.  For every nonarchimedean place $\lambda$ of $K_f$ dividing $\ell$, we denote by
$$
\rho = \rho_{f, \lambda} : \Gal (\oQ / F) \Map GL_2 ( K_{f, \lambda} )
$$
the associated semisimple, $K_f$-rational and integral Galois representation: 
See Deligne \cite{deligne-bourbaki}, Ohta \cite{ohta}, Carayol \cite{carayol}, Wiles \cite{wiles}, Taylor \cite{taylor-hilbert}, Blasius and Rogawski \cite{blasius-rogawski} and the references therein.
The Tate twist of its determinant $\det(\rho)(1)$ is a character of finite order. 

Let $G=G_{f, \lambda}$ be the Zariski closure of the image of $\rho_{f, \lambda}$ in $GL_2$ over $K_{f, \lambda}$.
Since we assume $\rho$ to be semisimple, the derived group of the connected component $(G^{\circ})^{\der} = [G^{\circ}, G^{\circ}]$
is a semisimple algebraic subgroup of $SL_2$, that is, either $SL_2$ or trivial.
If the reductive group $G^{\circ}_{f, \lambda}$ is a torus for some $\lambda$, we say that $f$ is of CM type.\footnote{The notion is independent of $\lambda$ and $\ell$ by a theorem of Serre,
cf. the argument in the proof of Theorem \ref{thm-cm}.}

The product
$$
\rho_{f, \ell} := \prod_{\lambda | \ell} \rho_{f, \lambda} : \Gal(\oQ / F) \Map \prod_{\lambda | \ell} GL_2( K_{f, \lambda} ) = (\Res^{K_f}_{\Q} GL_2) (\Ql)
$$
is $\Q$-rational and integral in the sense of Serre.  We denote by $G_{f, \ell}$ the Zariski closure over $\Ql$ of its image, and $G_{f,\ell}^{\circ}$ its connected component.

\subsection{Frobenius field and Ribet's argument}

\label{sec-ribet}

Following Ribet \cite{ribet}, for every finite extension $F'$ of $F$, we consider the Frobenius field:
$$
\Tr (\rho_{f, \lambda}, F'):= \Q \left( \left\{ \Tr \rho_{f, \lambda} ( \Frob_{\mathfrak{p}'} ) \right\}_{\mathfrak{p}'} \right) \le K_f
$$
where $\mathfrak{p}'$ ranges over the primes of $\mathcal{O}_{F'}$ coprime to $\disc(F)\mathfrak{n}\cdot\ell$; it is independent of $\lambda$, since the $\rho_{f, \lambda}$ form a strictly compatible system.  
Since $k_f = [K_f: \Q] < \infty$, there is the \textit{smallest Frobenius field} of $f$
$$
K_f^{\circ} \le K_f \,\,\,\,\mbox{ and } \,\,\,\, k_f^{\circ}:= [K_f^{\circ} :\Q].
$$
It is totally real, since $\det(\rho_{f, \lambda})(1)$ has finite order and the eigenvalues of $\rho_{f, \lambda}(\Frob_{\mathfrak{p'}})$ are Weil integers (see Lemma \ref{lem-Weil-number}).

We have thus the following algebraic groups:
$$
\xymatrix{
(\Res^{K_f^{\circ}}_{\Q} SL_2)_{\Ql} \ar@{^{(}->}[r] \ar@{^{(}->}[d] & ((\Res^{K_f^{\circ}}_{\Q} GL_2)^{\det\subseteq \Q^{\times}})_{\Ql} \ar@{^{(}->}[r] \ar@{^{(}->}[d]
    & (\Res^{K_f^{\circ}}_{\Q} GL_2)_{\Ql} \le GL_{2k_f^{\circ}, \Ql} \ar@{^{(}->}[d] \\
(\Res^{K_f}_{\Q} SL_2)_{\Ql} \ar@{^{(}->}[r] & ((\Res^{K_f}_{\Q} GL_2)^{\det\subseteq \Q^{\times}})_{\Ql} \ar@{^{(}->}[r]  & (\Res^{K_f}_{\Q} GL_2)_{\Ql} \le GL_{2k_f, \Ql} \\
}
$$
Here by the $\Q$-algebraic group $(\Res^{K_f}_{\Q} GL_2)^{\det\subseteq \Q^{\times}}$, we mean the following fibred product, which is often denoted by $G^{\ast}$ in the literature:
$$
\xymatrix{
G^{\ast} \ar[r] \ar[d] & \mathbb{G}_{m, \Q} \ar[d] \\
\Res^{K_f}_{\Q} GL_2 \ar[r]_-{\det} & \Res^{K_f}_{\Q} \mathbb{G}_{m, K_f},
}
$$
and similarly with $K_f$ replaced with $K_f^{\circ}$.

\begin{subproposition}[(Ribet)]
Suppose that $f$ is not of CM type.  Then
$$
G_{f, \ell}^{\circ} = ((\Res^{K_f^{\circ}}_{\Q} GL_2)^{\det\subseteq \Q^{\times}})_{\Ql}.
$$
\end{subproposition}
\begin{proof}
This amounts to showing that the (algebraic) Lie algebra $\mathfrak{g}$ of $G_{f, \ell}^{\circ}$ is equal to that of the right hand side.

The containment $\subseteq$ follows from the fact that, if $F'$ is any sufficiently large finite extension of $F$, 
then $\rho_{f, \lambda}(\Frob_{\mathfrak{p'}})$ has trace and determinant in $K_f^{\circ}$ for any prime $\mathfrak{p}'$ of $F'$ coprime to $\disc(F)\mathfrak{n}\ell$,
so that if $\lambda$ and $\lambda'$ are any $2$ primes of $K_f$ that lie over the same prime of $K_f^{\circ}$ over $(\ell)$, then $G^{\circ}_{f,\ell}$ is contained in the partial diagonal
of $(\Res^{K_f}_{\Q} GL_2)_{\Ql}$ where the $\lambda$- and the $\lambda'$-components are equal.

To prove the containment $\supseteq$, we first note that since $f$ is not of CM type, 
$\mathfrak{g}$ surjects onto each factor $\mathfrak{gl}_{2, \lambda^{\circ}}$, which contains $\mathfrak{sl}_{2, \lambda^{\circ}}$.
If $\lambda_1^{\circ}$ and $\lambda_2^{\circ}$ are distinct primes of $K_f^{\circ}$ lying over $\ell$, then the representations of $\mathfrak{g}$ in the $2$ factors are nonisomorphic,
since the representations of (germs of) $\Gal(\oQ/F)$ have different traces, and the image of $\mathfrak{g}$ in $\mathfrak{gl}_{2, \lambda_1^{\circ}} \times \mathfrak{gl}_{2, \lambda_2^{\circ}}$ contains $\mathfrak{sl}_{2} \times \mathfrak{sl}_2$.

Then by Goursat's Lemma (in the middle of the proof of Ribet \cite[Th. 4.4.10]{ribet}) the image of $\mathfrak{g}$ contains $\prod_{\lambda^{\circ}|\ell} \mathfrak{sl}_{2, \lambda^{\circ}}$,
where $\lambda^{\circ}$ ranges over the primes of $K_f^{\circ}$ lying over $(\ell)$.  Since the determinant on $G_{f, \lambda^{\circ}}^{\circ}$ is a dominant map onto $\Gm$, we get the desired equality.
\end{proof}


\begin{subdefinition}
Let $F^{\circ}$ be the Galois extension of $F$ cut out by two representations with finite image:
$$
\Gal (\oQ / F^{\circ}) = \ker ( \rho_{f, \ell}: \Gal(\oQ / F) \Map G_{f, \ell} (\Ql) / G^{\circ}_{f, \ell} (\Ql) )) \cap \ker ( \det(\rho_{f, \ell}) (1)).
$$
and let $\widetilde{F^{\circ}}$ be the compositum $F^{\circ}\widetilde{F}$.
\end{subdefinition}

\subsection{Variants of Sato-Tate equidistribution}

\label{sec-sst}

Let $i_1, \cdots, i_{k_f^{\circ}}$ denote the complete set of embeddings of $K_f^{\circ}$ into $\R$.
For each maximal ideal $\mathfrak{p}$ of $\mathcal{O}_{\widetilde{F^{\circ}}}$ coprime to $\disc(F)\cdot \mathfrak{n}\cdot \ell$, we let
$$
a_{\mathfrak{p}} = \Tr \rho_{f, \lambda} (\Frob_{\mathfrak{p}})
$$
and consider the set of vectors in $\R^{k_f^{\circ}}$:
$$
A(f)=\left\{ \left( \frac{i_1 (a_{\mathfrak{p}})}{\sqrt{\mathbb{N}\mathfrak{p}}}, \cdots, \frac{i_{k_f^{\circ}} (a_{\mathfrak{p}})}{\sqrt{\mathbb{N}\mathfrak{p}}} \right) \right\}_{\mathfrak{p}}.
$$

\begin{subdefinition}
We say that $f$ satisfies (SST) if $A(f)$ is equidistributed in the $k_f^{\circ}$-fold product of the Sato-Tate (half-circle) measure on $[-2, 2]$;
say that $f$ satisfies (RST) if $A(f)$ is equidistributed in a measure $\varphi(\mathbf{x}) d\mu_L(\mathbf{x})$, 
where $\varphi: [-2, 2]^{k_f^{\circ}} \Map \R_{\ge 0}$ is a continuous function and $d\mu_L$ is the Lebesgue measure.

For an integer $t\in [1, k_f^{\circ}]$, we say that $f$ satisfies ($t$-ST') if there exists a sequence $1\le j_1 < j_2 < \cdots < j_t \le k_f^{\circ}$ such that the projection
$\mathrm{pr}_{j_1, \cdots, j_t} (A(f)) \subseteq [-2, 2]^t$ is equidistributed in the $t$-fold product of the Sato-Tate measure on $[-2, 2]$;
we say that $f$ satisfies ($t$-ST) if for all sequences $\mathbf{j}$ of length $t$, $\mathrm{pr}_{\mathbf{j}}(A(f))$ is equidistributed in the $t$-fold product of the Sato-Tate measure.
\end{subdefinition}

We expect the strongest (SST) to be true; it fits into Serre's general framework of Sato-Tate equidistribution \cite[Chpt. 8]{serre-nxp}, for almost all primes $\ell$ (depending on $f$).
\footnote{The construction in \S 8.3 of \textit{op. cit.}, as stated, deals only with representations that come from the $\ell$-adic cohomology of algebraic varieties, 
but appears to use the condition only to the extent that they are rational and Hodge-Tate.

In case $d$ is odd or the automorphic representation corresponding to $f$ has a discrete series at some finite prime,
the Hodge-Tate condition (even the de Rham condition) for all $\ell$ follows from the motivic nature of the available constructions and theorems of Faltings \cite{faltings-jami},
see Blasius and Rogawski \cite{blasius-rogawski}.

In the general case, the Hodge-Tate condition is known for all but finitely many $\ell$: See Taylor \cite{taylor-hilbert-2}, where they are shown to be (even) crystalline.}
Namely, the compact Lie group attached to $\rho_{f, \ell}|_{\widetilde{F^{\circ}}}$ by Serre is the product of $k_f^{\circ}$ copies of $SU_2$
and the axioms $(A_1)$ and $(A_2)$ should hold.

When $t<k_f^{\circ}$, the mere conjunction of (RST) and ($t$-ST) does not imply (SST).

\begin{subremark}
The condition (SST) is stronger than the (usual) Sato-Tate equidistribution theorems available at the moment (see \cite{taylor-sato-tate-1}, \cite{taylor-sato-tate-2}, \cite{sato-tate-hilbert}).

In order to prove (SST) in the manner that the aforementioned results were obtained,
one would need to control the $L$-functions not only of the symmetric powers:
$$
\mathrm{Sym}^{m_j} (i_j(f))
$$
(through potential automorphy), but also of their tensor products:
$$
\mathrm{Sym}^{m_1} (i_1(f)) \otimes \cdots \otimes \mathrm{Sym}^{m_{k_f^{\circ}}} (i_{k_f^{\circ}} (f)).
$$
for all tuples $(m_1, \cdots, m_{k_f^{\circ}})$.  

The case of ($t$-ST), for $t\le 2$, looks accessible, see Harris \cite{harris}.
\end{subremark}

\subsection{Multisets and Newton Polygons}

\label{sec-semiring-newton}

Definition \ref{dfn-Pdki}, the operations $\otimes$, $\oplus$ and the partial order will be used in Theorem \ref{thm-main}.

We consider finite subsets with positive finite multiplicities, or simply \textit{multisets}, of $\Q$.  For example, $\{1/3, 2/3\}$ (each with multiplicity $1$) and $\{ 1/2, 1/2 \}$ (with multiplicity $2$).

\begin{subdefinition}
Let $S = \{ s_1, \cdots, s_m \}$ and $T = \{ t_1, \cdots, t_n \}$ be multisets.  Define the sum
$$
S\oplus T = \{ s_1, \cdots, s_m, t_1, \cdots, t_n \},
$$
the product
$$
S\otimes T = \{ s_i + t_j \}_{1\le i\le m, \,\, 1\le j\le n},
$$
and the dual
$$
S^{\vee} = \{ -s_1, \cdots, -s_m \}.
$$
Also, for $k> 0$, we write $S^{\oplus k} = S \oplus \cdots \oplus S$ and $S^{\otimes k} = S\otimes \cdots \otimes S$, repeated $k$ times.

The cardinality is denoted by $|S|$ or $\rank S$ (which is $m$ for the $S$ as above) and 
$$
\int S := s_1 + \cdots + s_m
$$ 
\end{subdefinition}

\begin{subproposition}
Multisets form a commutative semiring with involution, in which the empty set is the additive neutral element and $\{ 0\}$ the multiplicative identity element.

The map $S\mapsto |S|$ is a semiring homomorphism into the natural numbers and
$$
\int (S\oplus T) = \int S + \int T \,\,\, \mbox{ and } \,\,\, \int (S\otimes T) = |T| \int S + |S| \int T.
$$
\end{subproposition}

Given a multiset consisting of $a_1 \le \cdots \le a_n$ 
we form its Newton polygon emanating from $(0,0)$ with the slopes $a_1, \cdots, a_n$ (in this order).
Conversely, any finite Newton polygon emanating from the origin and with rational slopes uniquely determines a multiset of $\Q$.

From this point on, we will thus identify multisets with Newton polygons.  This allows us to impose a partial order on the class of multisets:
$$
S \le T \,\, \,\, \mbox{ if and only if } \,\, \,\, |S| = |T|\,\, \mbox{ and } \,\, \Newton(S) \le \Newton(T);
$$
the last meaning that $\Newton(T)$ lies on or above $\Newton(S)$.

\begin{subproposition}
Let $S \le S'$ and $T$ be three multisets. Then (1) $S \oplus T \le S' \oplus T$ and (2) $S \otimes T \le S' \otimes T$.  
If, in addition, $S$ and $S'$ end at the same point, then (3) $S^{\vee} \le S'^{\vee}$.
\label{prop-order}
\end{subproposition}
\begin{proof}

\noindent (1) By induction on $|T|$, we are reduced to the case where $T$ consists of $1$ element, say $T=\{ t\}$.
Twisting by $-t$ (i.e. taking $ \otimes \{ -t \}$) allows us to assume that $t=0$.
Enumerate $S$ and $S'$ in the  order:
$$
s_1 \le s_2 \le \cdots \le s_m \,\,\, \mbox{ and } \,\,\, s'_1 \le s'_2 \le \cdots \le s'_m
$$
and let $a$ and $b$ be such that:
$$
s_a < 0 \le s_{a+1}\,\,\, \mbox{ and } \,\,\, s'_b < 0 \le s'_{b+1};
$$
if all the $s_i$ are $\ge 0$ (resp. $<0$), then we let $a:= 0$ (resp. $a:=m$), and similarly for $b$.

If we define $\Sigma_S(i):= s_1 + \cdots + s_i$ for $i\in [0, m]$, the condition $S\le S'$ becomes
$$
\Sigma_S(i) \le \Sigma_{S'}(i) \,\,\, \mbox{ for all } \,\, i\in [0, m].
$$
We need to prove 
\begin{equation}
\label{eqn-newton-sum}
\Sigma_{S\oplus \{0\}} (i) \le \Sigma_{S'\oplus \{ 0 \}} (i) \,\,\, \mbox{ for all } \,\, i\in [0, m+1].
\end{equation}

\noindent (1a) Suppose that $a<b$.  Then for $i\in [0, a] \cup [b+1, m+1]$, (\ref{eqn-newton-sum}) is clearly satisfied.  For $i\in [a+1, b]$,
since $0$ is inserted into $S'$ at the $(b+1)$-th place, we have 
\begin{eqnarray*}
\Sigma_{S'\oplus\{0\}}(i) & = & \Sigma_{S'\oplus\{0\}}(b+1) - ( s'_{i+1} + \cdots + s'_b + 0 ) \\
 & = & \Sigma_{S'}(b) - (s'_{i+1} + \cdots + s'_b) \ge \Sigma_{S'}(b)
\end{eqnarray*}
while, since $0$ is inserted into $S$ at the $(a+1)$-th place, we have
\begin{eqnarray*}
 \Sigma_{S \oplus \{ 0\} } (i) & = & \Sigma_{S}(i-1) = \Sigma_S(b) - ( s_{i} + \cdots + s_b ) \le \Sigma_S(b).
\end{eqnarray*}
and we get (\ref{eqn-newton-sum}).

\noindent (1b) Suppose that $a\ge b$.  Then again for $i\in [0,b] \cup [a+1, m+1]$, (\ref{eqn-newton-sum}) is trivially satisfied.  For $i\in [b+1, a]$, we have this time:
\begin{eqnarray*}
\Sigma_{S'\oplus\{0\}} (i) & = & \Sigma_{S'}(i-1) = \Sigma_{S'}(b) + (s'_{b+1} + \cdots + s'_{i-1}) \ge \Sigma_{S'}(b) \\
\Sigma_{S\oplus\{0\}} (i) & = & \Sigma_S(i) = \Sigma_S(b) + (s_{b+1} + \cdots + s_i) \le \Sigma_{S'} (b)
\end{eqnarray*}
This completes the proof of (1).

\vspace{12 pt}

\noindent (2) By decomposing $T$ into singletons and using the distributive law, we deduce (2) from (1).

\vspace{12 pt}

\noindent (3) The duals $S^{\vee}$ and $S'^{\vee}$ are enumerated:
$$
-s_m \le -s_{m-1} \le \cdots \le -s_1 \,\,\, \mbox{ and } \,\,\, -s'_m \le -s'_{m-1} \le \cdots \le -s'_1.
$$
The assumption that $S$ and $S'$ end at the same point means that $\Sigma_S(m)=\Sigma_{S'}(m)$.  Thus
$$
\Sigma_{S^{\vee}}(i) =  \Sigma_S(m-i) - \Sigma_S(m) \le \Sigma_{S'}(m-i) - \Sigma_{S'}(m) = \Sigma_{S'^{\vee}}(i)
$$
for all $i\in [1, m]$, and this completes the proof of the Proposition.
\end{proof}

\begin{subremark}
In view of (3), one may want to consider the more restrictive partial order:
$$
S \le' T \,\,\, \mbox{ if and only if } \,\,\, |S| = |T|, \int S = \int T, \mbox{ and } \Newton(S) \le \Newton(T),
$$
so as to make the involution $S\mapsto S^{\vee}$ order-preserving.  Below, we use the partial order $\le$ only in the case where $\le'$ also applies.
\end{subremark}

We are particularly interested in:
\begin{subdefinition}
By the partially ordered semiring of \textbf{integral Newton polygons}, we mean the subsemiring of multisets whose Newton polygons have integral breaking points.
\end{subdefinition}
The following polygons appear in the statement of Theorem \ref{thm-main}.
\begin{subdefinition}
\label{dfn-Pdki}
Let $d\ge 1$, $k\ge 1$ and $i\in [0, k]$ be integers.  We define the multiset (and the corresponding Newton polygon):
$$
P(d; k, i) := \left( \{ 0, 1 \} ^{\otimes d} \right)^{\oplus (k-i)} \oplus \left( \{ 1/2, 1/2 \} ^{\otimes d} \right)^{\oplus i}
$$
\end{subdefinition}
The Newton polygon of $P(d; k,i)$ has integral breaking points.  By Proposition \ref{prop-order}, we have
$$
P(d; k, i) \le P(d; k, j) \,\,\,\, \mbox{ if and only if } \,\,\,\, i \le j.
$$

\subsection{Interaction of $F$ and $K$: Slope and Bisection}
\label{sec-interplay}

Let $G$ be a group acting on a finite set $X$.

\begin{subdefinition}
 By the \textbf{length of maximal parts} of $g\in G$ on $X$, which we denote by $\lambda(g, X)$, we mean the largest of the cardinalities of the $g$-orbits in $X$.
We define $\lambda(G, X)$ as the supremum of $\lambda(g, X)$, as $g$ ranges over $G$.
\end{subdefinition}

\begin{subdefinition}
We say that an element $g\in G$ \textbf{bisects} 
$X$ if $g$ has exactly $2$ orbits in $X$ and the orbits have the same number of elements.
\end{subdefinition}

Let $F$ be a (ground) number field, $K$ a (coefficient) number field, and $F^s$ an algebraic closure of $F$.
The Galois group $G:= \Gal(F^s/F)$ acts continuously on the discrete set 
$$
X:= \Hom ( K, F^s )
$$
of field embeddings of $K$ into $F^s$.  

\begin{subdefinition}
We define
$$
\lambda_{F}(K) := \lambda(\Gal(F^s/F), \Hom(K, F^s)).
$$
When $F=\Q$, we drop $F$ from the notation and write $\lambda(K)$.
\end{subdefinition}

In more concrete terms: When $F=\Q$, the Galois group of the normal closure $\widetilde{K}$ of $K/\Q$ determines $\lambda(K)$.
For example, if the group is the full symmetric group of degree $[K:\Q]$ or the cyclic group of $[K:\Q]$ elements, then $\lambda(K) = [K:\Q]$.
If the group is the alternating group, then $\lambda(K) = [K:\Q]$ (resp. $=[K:\Q]-1$) if $[K:\Q]$ is odd (resp. even).

The notion of bisection is also determined by the Galois group; for example, the alternating group of even degree and the Klein $4$-group acting on itself by translation contain bisecting elements.

For general $F$, one needs to look at the action of the subgroup
$\Gal ( \oQ / F ) \le \Gal( \oQ / \Q).
$

\begin{subdefinition}
\label{def-slope}
Given two number fields  $F$ and $K$, we define 
the \textbf{slope} of $K$ over $F$:
$$
\sigma_{F} (K) = 1- \frac{\lambda_{F} (K)}{[K: \Q]}\,\,\,\, \in [0, 1) \cap \Q.
$$
When $F=\Q$, we write $\sigma(K)=\sigma_{\Q}(K)$.
\end{subdefinition}
We call $\sigma$ the slope, in view of the following `semistability' property, formally analagous to that of Harder and Narasimhan for vector bundles on curves, in the variable $K$:

\begin{subproposition}
\label{prop-semistability}
Let $K'$ be a subfield of $K$ and let $n=[K:K']$.  Then
$$
n \cdot \lambda (g, \Hom(K', F^s) ) \ge \lambda (g, \Hom(K, F^s))
$$
for all $g\in G = \Gal(F^s/F)$ and therefore
$$
\sigma_F(K') \le \sigma_F(K).
$$
\end{subproposition}
\begin{proof}
Let $X = \Hom (K, F^s)$ and $X':= \Hom(K', F^s)$.  Then we have a surjective map of $G$-sets:
$$
X \Map X'
$$
obtained by restriction.  Since $K/K'$ is separable, each fibre has exactly $n$ elements. 

The first inequality follows from inspecting the images of the $g$-orbits in $X$, and the second follows from the first by definition.
\end{proof}

In the variable $F$, we trivially have
$$
\sigma_{F}(K) \le \sigma_{F'} (K) \,\,\,\, \mbox{ if } \,\,\,\, F \subseteq F'.
$$

The following is useful in computing $\sigma$ and finding bisecting elements in practice.
\begin{subproposition}
\label{prop-interact}
Let $F$ and $K$ be two number fields, with the respective normal closures $\widetilde{F}$ and $\widetilde{K}$ over $\Q$.
\begin{enumerate}
 \item[(1)] If $[K:\Q]$ is a prime number not dividing $[F:\Q]$, then $\sigma_F(K) = 0$.
 \item[(2)] If $\Gal(\widetilde{K}/\Q)$ is the symmetric group of degree $[K:\Q]$ and if $[\widetilde{F}:\Q]$ is odd, then $\sigma_{\widetilde{F}}(K)= 0$.
 \item[(3)] Suppose that $\widetilde{K}$ is linearly disjoint from $\widetilde{F}$ over $\Q$ (which is the case, for example, if $F$ and $K$ have coprime discriminants). 
 Then $\sigma_F(K) = \sigma_{\Q} (K)$, and $\Gal(\oQ/F)$ possesses an element bisecting $\Hom(K, \oQ)$ exactly when $\Gal(\oQ/\Q)$ does.
 \end{enumerate}
\end{subproposition}

\begin{proof}  (1) Recall that $\Gal(\widetilde{K}/\Q)$ acts transitively on $\Hom(K, \Q)$, and therefore has order divisible by $p=[K:\Q]$.
Since the image of $\Gal(\oQ/F)$ in $\Gal(\widetilde{K}/\Q)$ (via $\Gal(\oQ/\Q)$) has index dividing $[F:\Q]$, the image also has order divisible by $p$.
Therefore the image contains an element $g$ of exact order $p$, which has $\lambda(g, \Hom(K, \Q))=p=[K:\Q]$.

(2) Use the fact that a symmetric group has no proper normal subgroup of odd index to deduce that the image of $\Gal(\oQ/\widetilde{F})$ in $\Gal(\widetilde{K}/\Q)$ is the full symmetric group.

(3) By assumption, the natural map
$$
(\phi_1, \phi_2): \Gal(\oQ / \Q) \Map \Gal(\widetilde{F}/\Q) \times \Gal(\widetilde{K}/\Q)
$$
is surjective and by definition $\Gal(\oQ / F) \supseteq \Gal(\oQ/\widetilde{F}) = \ker(\phi_1)$.  Therefore $\Gal(\oQ/F)$ surjects onto $\Gal(\widetilde{K}/\Q)$.
The statements about the slope and bisecting elements follow from this.
\end{proof}

\subsection{Zariski density and Haar density (Serre)}

We will use the following in the proof of Theorem \ref{thm-main}.

\begin{sublemma}
\label{lem-serre}
Let $E$ be a finite extension of $\Q_{\ell}$, $G$ a connected algebraic group over $E$, $\Gamma \le G(E)$ a compact and Zariski dense subgroup, and
$$
\varphi: G \Map \mathbb{A}^1_E
$$
a regular morphism of algebraic varieties that is constant on the conjugacy classes.

Then for any finite subset $S$ of $E=\mathbb{A}^1(E)$ that does not contain $\varphi(1_G)$, the subset $\Gamma \cap \varphi^{-1}(S)$ has Haar measure $0$ in $\Gamma$.

In particular, if $\Gamma$ is the image of a continuous representation $\rho$ of $\Gal(\oQ/ F)$ unramified outside a finite set of primes of the number field $F$, 
then the set of primes $\mathfrak{p}$ in $F$ such that $\varphi(\rho(\Frob_{\mathfrak{p}})) \in S$ has (natural) density $0$.
\end{sublemma}
\begin{proof}
This follows from Serre \cite[Prop.5.12]{serre-nxp} applied to $Z:= \varphi^{-1}(S)$, which is a proper algebraic subset of $G$ by the assumptions and has Zariski density $0$ by definition.
\end{proof}

\section{Main Theorems}

\subsection{Non CM case}

The notation and terminology in the following theorem are explained in the previous preparatory section. References to the precise subsections are provided as they arise.

\begin{subtheorem}
\label{thm-main}
Let $f$ be a new normalised Hilbert eigencuspform of level $\mathfrak{n}\subseteq \mathcal{O}_F$ and parallel weight $(2, \cdots, 2)$, and suppose that it is not of CM type (\S \ref{sec-notation}).

Denote by $M(f)$ the Andr\'e motive (see Proposition \ref{rmk-hecke-andre-motives}), whose realisations give the part of the intersection cohomology of the Hilbert modular variety
corresponding to $\{\sigma(f)\}_{\sigma}$, where $\sigma$ ranges over all the embeddings of $K_f$ into $\oQ$.

\begin{enumerate}
 \item[(1)] (Analogue of the Katz Conjecture) For all rational primes $p$ coprime to $\disc(F)\cdot \mathfrak{n}\cdot \ell$, we have
 $$
 \Newton ( \Frob_p |_{M(f)} ) \ge \HodgeP (M(f)).
 $$
 Moreover, if $p$ splits completely in $F$ (equivalently in $\widetilde{F}$) and $p$ is unramified in $K_f$, then there exists an integer $k(p)\in[0, k_f]$ such that
 $$
 \Newton ( \Frob_p |_{M(f)}) = P(d; k_f, k(p)).
 $$
 (Here $k_f = [K_f:\Q]$ and we refer to Definition \ref{dfn-Pdki} for the right hand side.)
 
 In the remaining parts, we only consider the primes splitting completely in $F$ and unramified in $K_f$.
 \item[(2)] For a principally abundant set of primes $p$, we have
 $$
 k(p) \le \frac{k_f}{2}.
 $$
 \item[(3)] If $k_f^{\circ} = [K_f^{\circ}: \Q] \le 2$ ($K_f^{\circ}$ is defined in \S \ref{sec-ribet}), then for a principally abundant set of primes $p$, we have $k(p)=0$,
 that is, the Newton and Hodge-Tate polygons coincide.
 \item[(4)] For an abundant set of primes $p$, we have ($\sigma$ defined in \S \ref{sec-interplay})
 $$
 k(p) \le k_f \cdot \mathrm{min} \left\{ 1/2 , \sigma_{\widetilde{F}}(K_f) \right\}
 $$
 \item[(4bis)] For an abundant set of primes $p$, we have ($\widetilde{F^{\circ}}$ defined in \S \ref{sec-ribet})
 $$
 k(p) \le k_f \cdot \mathrm{min} \left\{ 1/2 , \sigma_{\widetilde{F^{\circ}}}(K_f^{\circ}) \right\}
 $$
 \item[(5)] If $k_f^{\circ}$ is even, suppose that $f$ satisfies (RST) (resp. ($t$-ST' for an integer $t\ge 1$), as defined in \S \ref{sec-sst}. 
 Then for a principally abundant set of primes $p$ (resp. for an abundant set of primes $p$), we have
 $$
 k(p) \le \frac{k_f}{k_f^{\circ}} \left\lfloor \frac{k_f^{\circ} -1}{2} \right\rfloor;
 $$
 in particular, $k(p) < k_f / 2$.
 \item[(6)] Suppose that $f$ satisfies (RST) and that an element of $\Gal(\oQ / \widetilde{F^{\circ}})$ 
 bisects  $\Hom(K_f^{\circ}, \oQ)$ (\S \ref{sec-interplay}).
 Then for an abundant set of primes $p$, we have $k(p)=0$, i.e., the Newton and Hodge-Tate polygons coincide.
\end{enumerate}

\end{subtheorem}

\begin{proof}
We have already fixed a rational prime $\ell$ that splits completely in $K_f$.  Now for all rational primes $p \neq \ell$, we fix once and for all an isomorphism
$
\xymatrix{
\oQ_{\ell} \ar[r]^-{\sim} & \oQ_p,
}
$
and pull back the $p$-adic valuation on the target to get a rank-$1$ (discontinuous) valuation $v_p$ on $\oQ_{\ell}$, normalised by $v_p(p) = 1$.

To prove (1), we may pass to $\widetilde{F}$, and consider $\Frob_{\mathfrak{p}}$ for any prime $\mathfrak{p}$ lying over $p$,
because doing so does not change the Newton polygon.

The key fact (from Brylinski and Labesse \cite{BL}) that we use from the constructions 
(see Deligne \cite{deligne-bourbaki}, Ohta \cite{ohta}, Carayol \cite{carayol}, Wiles \cite{wiles}, Taylor \cite{taylor-hilbert} and Blasius and Rogawski \cite{blasius-rogawski} and the references therein)
of the Galois representations associated with the $\{ \sigma(f) \}$ is the following:
the $\ell$-adic \'etale realisation $M(f)_{\ell}$ of $M(f)$ is the direct sum of the tensor inductions (see Curtis and Reiner \cite[\S 80C]{curtis-reiner}):
\begin{equation}
\label{eqn-tensor-ind-gen}
\bigoplus_{\sigma\in \Hom(K_f, \Ql)} \tensorInd^{\Gal(\oQ / \Q)}_{\Gal(\oQ / F)} (\rho_{\sigma(f), \lambda}).
\end{equation}
This implies that for transversals (coset representatives) $g_1, \cdots, g_d \in \Gal(\oQ/\Q)$ of $\Gal(\oQ/\Q)$ modulo $\Gal(\oQ / F)$, we have:
For any element $\gamma$ in the finite index normal subgroup $\Gal(\oQ / \widetilde{F})$ of $\Gal(\oQ / F)$, the action of $\gamma$ on $M(f)_{\ell}$ is given by
\begin{equation}
\label{eqn-tensor-ind-explicit}
\bigoplus_{\sigma} \left( \rho_{\sigma(f), \lambda} (g_1 \gamma g_1^{-1}) \otimes \rho_{\sigma(f), \lambda} ( g_2 \gamma g_2^{-1} ) \otimes \cdots \otimes \rho_{\sigma(f), \lambda} ( g_d \gamma g_d^{-1} ) \right)
\end{equation}
Let $K':= K_f (\Frob_{\mathfrak{p}})$ be the splitting field over $K_f$ of the polynomial
\begin{equation}
\label{eqn-quadratic}
X^2 - \Tr ( \rho(\Frob_{\mathfrak{p}})) X + \det (\rho(\Frob_{\mathfrak{p}})),
\end{equation}
and let 
$
R_{\mathfrak{p}} = \{ \alpha_{\mathfrak{p}}, \beta_{\mathfrak{p}}\} \subset K'
$ 
be the roots. Each embedding $\sigma: K_f \Map \Ql$ extends to (at most $2$) embeddings $\sigma':K' \Map \oQ_{\ell}$,
and the image $\sigma'(R_{\mathfrak{p}}) \subset \oQ_{\ell}$ is independent of the choice $\sigma'$.
Therefore we can unambiguously form the multiset of slopes\footnote{this may not have integral breaking points}:
$$
S_{\mathfrak{p}, \sigma} = \left\{ \frac{v_p(\sigma'(\alpha_{\mathfrak{p}}))}{v_p(\mathbb{N}\mathfrak{p})}, \frac{v_p(\sigma'(\beta_{\mathfrak{p}}))}{v_p(\mathbb{N}\mathfrak{p})} \right\}.
$$
Since $\Tr(\rho(\Frob_{\mathfrak{p}}))$ is an algebraic integer, its $p$-adic valuation is $\ge 0$.
Also, $\det(\rho(\Frob_{\mathfrak{p}}))$ is $\mathbb{N}\mathfrak{p}$ times a root of unity, so its $p$-adic valuation is equal to $v_p(\mathbb{N}\mathfrak{p})$.
These two facts imply the inequalities on the $\sigma$-slopes of (\ref{eqn-quadratic}) for all $\sigma$:
\begin{equation}
\label{eqn-2-bounds}
\{ 0, 1 \} \le S_{\mathfrak{p}, \sigma} \le \{ 1/2, 1/2 \}
\end{equation}
(see \S \ref{sec-semiring-newton} for the partial order by Newton polygon). Moreover, if $\mathbb{N}\mathfrak{p}=p$ (in particular, if $p$ splits completely in $F$)
and $p$ is unramified in $K_f$, then one of the two inequalities must be an equality.

In view of the description of cohomology in terms of tensor induction (\ref{eqn-tensor-ind-gen}) and (\ref{eqn-tensor-ind-explicit}), we have
$$
\Newton(\Frob_{\mathfrak{p}|_{M(f)}}) = \bigoplus_{\sigma\in \Hom(K_f, \Ql)} S_{\mathfrak{p}, \sigma}^{\otimes d},
$$
and since 
$$
\HodgeP (M(f)) = HP(M(f)) = P(d; k_f, 0) = (\{0, 1\}^{\otimes d})^{\oplus k_f},
$$
we have proven (1).

\begin{subremark}
We also get a bound for the denominators of the slopes: they are divisors of integers in the interval $[1, \widetilde{d}]$, or equal to $2$.
\end{subremark}

In order to proceed further, we first recall the following known fact (generalised Ramanujan-Petersson conjecture, see Taylor \cite{taylor-icm1994} and Blasius\cite{blasius-ramanujan}):
\begin{sublemma}
\label{lem-Weil-number}
The roots $\alpha_{\mathfrak{p}}$ and $\beta_{\mathfrak{p}}$ are $\mathbb{N}\mathfrak{p}$-Weil integers of weight $1$.
\end{sublemma}
\begin{proofarg}{of Lemma}
In case $d$ is odd or the automorphic representation $\pi_f$ corresponding to $f$ is a discrete series representation at some finite prime, this follows from
the essentially motivic nature of some of the constructions, see Blasius and Rogawski \cite{blasius-rogawski}, together with Deligne's proof of the Weil conjectures \cite{deligne-weil}.

This can be proved in the general case, and only with the a priori non motivic constructions of Wiles \cite{wiles} and Taylor \cite{taylor-hilbert}.  
Note that by the description (\ref{eqn-tensor-ind-gen}), the algebraic integers
$
\sigma'(\alpha_{\mathfrak{p}}^d) \, \mbox{ and } \, \sigma'(\beta_{\mathfrak{p}}^{d})
$
for any embedding $\sigma'$ of $K'$ into $\oQ_{\ell}$, are eigenvalues of $\Frob_{\mathfrak{p}}$ acting on the $IH^d$ of the Baily-Borel compactification $X^{BB}(\mathfrak{n})$ of the Hilbert modular variety.

Now this last variety admits a surjective, generically finite map from a projective smooth toroidal compactification over $\Z_{(p)}$ \cite{rapoport}. 

By the decomposition theorem for perverse sheaves \cite{bbd}, the $2$ algebraic integers therefore appear as eigenvalues of $\Frob_{\mathfrak{p}}$ in the $H^d$ of the projective smooth variety.
Then by Deligne's proof of the Weil conjectures, they have all the archimedean absolute values $= (\mathbb{N}\mathfrak{p})^{d/2}$.  By taking the $d$-th root, we get the Lemma.
\end{proofarg}

\vspace{12 pt}

From this point on, we assume that $\mathfrak{p}$ is a prime of absolute degree $1$ over $(p)$ (in addition to being coprime to $\disc(F)\cdot\mathfrak{n}\cdot\ell$).
We also assume that $p$ is unramified in $K_f$.

\noindent (2) Now we look more closely at
$$
a_{\mathfrak{p}} := \Tr ( \rho_{f, \lambda}(\Frob_{\mathfrak{p}})).
$$
By the assumption that $f$ is not of CM type, the image of $\Gal(\Q / \widetilde{F})$ under $\rho_{f, \lambda}$ is Zariski dense in $GL_2$ over $K_{f, \lambda}$.
Since $\Tr$ is a regular morphism of the algebraic variety $GL_2$ into $\mathbb{A}^1$ and takes value $2\neq 0$ at the identity element $I_2$,
the set of primes $\mathfrak{p}$ of $\widetilde{F}$ such that $a_{\mathfrak{p}} = 0$ has density zero by Lemma \ref{lem-serre}.  We exclude them from this point on.

Let $\wp_1, \cdots, \wp_m$ be the primes of $K_f$ lying over $(p) = \mathfrak{p} \cap \Z$, and write the ideal factorisation 
\begin{equation}
\label{eqn-factorisation}
a_{\mathfrak{p}} \cdot \mathcal{O}_{K_f} = \wp_1^{e_1} \cdots \wp_k^{e_m} \cdot I,
\end{equation}
where $I$ is an integral ideal coprime to $p$, and we carry on with the argument preceeding the Lemma.

For each embedding $\sigma: K_f \Map \Ql$, let $\wp_{i(\sigma)}$ be the inverse image in $\mathcal{O}_{K_f}$ of the maximal ideal of the integral closure $\overline{\Z_p}\subset \oQ_p$ of $\Z_p$
under the composite of $\sigma$ with the fixed $\oQ_{\ell}\simeq \oQ_p$:
$$
\xymatrix{
\mathcal{O}_{K_f}\subset K_f \ar[r]^-{\sigma} & \oQ_{\ell} \simeq \oQ_p \supset \mathfrak{m}_{\overline{\Z_p}}.
}
$$
Then $S_{\mathfrak{p}, \sigma}$ is the Newton polygon of the polynomial (obtained by applying $\sigma$ to (\ref{eqn-quadratic})) with respect to $v_p$ on $\Ql$:
$$
X^2 - \sigma(a_{\mathfrak{p}}) X + \sigma (\det(\rho(\Frob_{\mathfrak{p}})))
$$
and as such is equal to:
$$
S_{\mathfrak{p}, \sigma} = \begin{cases}
                            \{ 0, 1 \} & \mbox{ if } e_{i(\sigma)} = 0 \\
                            \{ 1/2, 1/2 \} & \mbox{ if } e_{i(\sigma)}>0
                           \end{cases}
$$
Let $f_1, \cdots, f_m$ be the degrees of the residue class extensions:
$$
f_i := \dim_{\Fp} \mathcal{O}_{K_f} / \wp_i.
$$
Then we have
$$
\Newton (\Frob_{\mathfrak{p}}|_{M(f)}) = P(d; k_f, k(p)),
$$
where $k(p)$ is the sum of those $f_i$ for which $e_i>0$.

Since $a_{\mathfrak{p}} \neq 0$, we may apply the product formula.  By the Lemma, we have
$$
\prod_{v|\infty} \| a_{\mathfrak{p}} \|_v \le (2\sqrt{p})^{k_f},
$$
while by the factorisation (\ref{eqn-factorisation}):
$$
\prod_{v|p} \| a_{\mathfrak{p}} \|_v = p^{-\sum_{i=1}^m e_i f_i}
$$
and
$$
\prod_{v\nmid p, \infty} \|a_{\mathfrak{p}} \|_v = \prod_{v\nmid p, \infty} \| I \|_v \le 1.
$$
Therefore
$$
1 = \prod_{v} \| a_{\mathfrak{p}} \| \le 2^{k_f} p^{\frac{k_f}{2} - \sum_i e_i f_i},
$$
which implies for all $p>2^{2k_f}$
$$
\frac{k_f}{2} \ge \sum_{i=1}^m e_i f_i \ge \sum_{i: e_i>0} f_i = k(p).
$$
This proves (2).

\vspace{12 pt}

For (3), we assume in addition that $p$ splits completely in $\widetilde{F^{\circ}}$, and choose a prime $\mathfrak{p}$ of $\widetilde{F^{\circ}}$ lying over $p$,
so that by definition 
$$
\rho_{f, \ell}(\Frob_{\mathfrak{p}}) \in G_{f,\ell}^{\circ}(\Ql) \,\,\, \mbox{ and } \,\,\, a_{\mathfrak{p}} \in \mathcal{O}_{K_f^{\circ}}.
$$

If $k_f^{\circ}=1$, that is, if $K_f^{\circ} = \Q$, then $a_{\mathfrak{p}}\in \Z$ and $|a_{\mathfrak{p}}|_{\infty} < 2\sqrt{p}$.
As soon as $p\ge 5$, the only way $p| a_{\mathfrak{p}}$ is then $a_{\mathfrak{p}}=0$, which we have excluded above.

Suppose therefore that $k_f^{\circ} = 2$, and consider the homomorphisms
$$
\xymatrix{
\Gal ( \Q / \widetilde{F^{\circ}}) \ar[r]^-{\rho^{\circ}} & \left( (\Res^{K_f^{\circ}}_{\Q} GL_2)^{\det \subseteq \Q^{\times}}\right) (\Ql)
  \ar@{^{(}->}[r]^-{\iota} \ar[d]^-{\det} & GL_4(\Ql) \\
 & \Ql^{\times} &
}
$$
and the regular map of algebraic varieties
\begin{equation}
\Tr( \wedge^2(\iota) \otimes {\det}^{-1} ): ((\Res^{K_f^{\circ}}_{\Q} GL_2)^{\det \subseteq \Q^{\times}}) \otimes \Ql \Map \mathbb{A}^1_{\Ql},
\label{eqn-reg-map-katz}
\end{equation}
where $\wedge^2(\iota)$ takes values in $GL_6$ and $\det$ takes values in $\Gm$.

In order to prove (3), we find a set of primes $\mathfrak{p}$ of $\widetilde{F^{\circ}}$ of density $=1$ such that $a_{\mathfrak{p}}$ is not divisible by any prime of $K_f^{\circ}$ lying over $(p)$.
If $p$ is inert in $K_f^{\circ}$, the bound (2) suffices, and we exclude the finitely many primes that are ramified in $K_f^{\circ}$, so we assume that $p$ splits:
\begin{equation}
\label{eqn-p-splits}
p\cdot \mathcal{O}_{K_f^{\circ}} = \wp_1 \wp_2, \,\,\,\, \mbox{ where } \wp_1 \neq \wp_2 \mbox{ are primes.}
\end{equation}

By (2), we may assume that at most one of the $2$ primes can divide $a_{\mathfrak{p}}$, say $a_{\mathfrak{p}} \in \wp_1$ but $a_{\mathfrak{p}} \not\in \wp_2$.
Let $\epsilon$ be the nontrivial field automorphism of $K_f^{\circ}$, and let $\alpha_{\mathfrak{p}}$ be a root of the polynomial
$
X^2 - a_{\mathfrak{p}} X + p = 0
$
and let $\alpha'_{\mathfrak{p}}$ be a root of
$
X^2 - \epsilon(a_{\mathfrak{p}} ) X + p = 0.
$

Then the $4$ eigenvalues of $(\iota \circ{\rho^{\circ}}) (\Frob_{\mathfrak{p}})$ are 
$\left\{ \alpha_{\mathfrak{p}}, \frac{p}{\alpha_{\mathfrak{p}}}, \alpha'_{\mathfrak{p}}, \frac{p}{\alpha_{\mathfrak{p}}} \right\}$,
and the $6$ eigenvalues of $(\wedge^2 \iota \circ{\rho^{\circ}}) (\Frob_{\mathfrak{p}})$ are
$
\left\{ p, \, \alpha_{\mathfrak{p}} \alpha'_{\mathfrak{p}}, \, \alpha_{\mathfrak{p}}\frac{p}{\alpha'_{\mathfrak{p}}}, \, \frac{p}{\alpha_{\mathfrak{p}}}\alpha'_{\mathfrak{p}}, \,
  \frac{p}{\alpha_{\mathfrak{p}}} \frac{p}{\alpha'_{\mathfrak{p}}}, \, p \right\}.
$
Therefore the value of (\ref{eqn-reg-map-katz}) at $\Frob_{\mathfrak{p}}$ is
$$
\frac{1}{p} \left( 2p + (\alpha_{\mathfrak{p}} + \frac{p}{\alpha_{\mathfrak{p}}}) (\alpha'_{\mathfrak{p}} + \frac{p}{\alpha'_{\mathfrak{p}}}) \right) 
    = 2 + \frac{a_{\mathfrak{p}}\epsilon(a_{\mathfrak{p}})}{p} \in \Z
$$
Now by the Lemma \ref{lem-Weil-number}, which implies
$$
|a_{\mathfrak{p}}|_{\infty}, |\epsilon(a_{\mathfrak{p}})|_{\infty} \le 2\sqrt{p}
$$
and by our assumption that $p$ is unramified in $K_f^{\circ}$, which implies that the inequalities are strict, we have
\begin{equation}
\Tr( \wedge^2(\iota) \otimes {\det}^{-1} ) (\Frob_{\mathfrak{p}}) \in [-1,5] \cap \Z.
\label{eqn-katz-trace}
\end{equation}

Therefore $a_{\mathfrak{p}}$ does not belong to any prime of $K_f^{\circ}$ lying over $(p)$, as soon as we avoid (\ref{eqn-katz-trace}).
But since
$$
\Tr( \wedge^2(\iota) \otimes {\det}^{-1} ) (1) = 6
$$
(here $1$ denotes the identity element in the group $G^{\circ}_{f, \ell}$), the set of $\mathfrak{p}$ for which (\ref{eqn-katz-trace}) holds has density $0$ by Lemma \ref{lem-serre}.
Now if $\mathfrak{p}$ avoids (\ref{eqn-katz-trace}), then
$$
a_{\mathfrak{p}} \mathcal{O}_{K_f} = (a_{\mathfrak{p}} \mathcal{O}_{K_f^{\circ}}) \mathcal{O}_{K_f}
$$
is coprime to $(p)$, and we have $k(p) = 0$.  This completes the proof of (3).

\vspace{12 pt}

For the sake of continuity in exposition, we treat the conditional (5) before the unconditional (4).  By an argument similar to that for (2), but applied to
the restriction $\rho^{\circ}: \Gal ( \oQ / \widetilde{F^{\circ}} ) \Map GL_2( K^{\circ}_{f, \lambda^{\circ}})$ (where $\lambda^{\circ} = \lambda \cap K^{\circ}_f$),
for a prime $\mathfrak{p}$ of density $1$ in $\widetilde{F^{\circ}}$, if we write
$$
a_{\mathfrak{p}} \mathcal{O}_{K^{\circ}_f} = \wp_1^{e_1} \cdots \wp_{m'}^{e_{m'}} I
$$
where $\wp_1, \cdots, \wp_{m'}$ are the primes of $K^{\circ}_f$ lying over $p$ and $I$ is coprime to $p$, then we have
\begin{equation}
\label{eqn-proof-5}
\sum_{i: e_i>0} \dim_{\Fp} ( \mathcal{O}_{K^{\circ}_f} / \wp_i ) \le \frac{k_f^{\circ}}{2},
\end{equation}
and
$$
k(p) = \frac{k_f}{k_f^{\circ}} \sum_{i: e_i>0} \dim_{\Fp} ( \mathcal{O}_{K^{\circ}_f} / \wp_i ).
$$

If $k_f^{\circ}$ is odd, then (\ref{eqn-proof-5}) trivially implies (5), so we assume that $k_f^{\circ}$ is even.

Now if the equality holds in (\ref{eqn-proof-5}), then we necessarily have $e_i=1$ whenever $e_i>0$, and by Lemma \ref{lem-Weil-number} and the product formula
$$
(2\sqrt{p})^{k_f^{\circ}} \ge \prod_{v|\infty} \| a_{\mathfrak{p}} \|_v = \left( \prod_{v|p} \|a_{\mathfrak{p}}\|_v \cdot \prod_{v\nmid p\infty} \|a_{\mathfrak{p}}\|_v\right)^{-1} \in p^{k_f^{\circ}/2}\Z.
$$
In other words, if $i_1, \cdots, i_{k_f^{\circ}}$ are the real embeddings of $K_f^{\circ}$, then we have
$$
\prod_{j=1}^{k_f^{\circ}} \frac{i_j (a_{\mathfrak{p}})}{\sqrt{p}} \in \Z \cap [-2^{k_f^{\circ}}, 2^{k_f^{\circ}} ].
$$
In $\R^{k_f^{\circ}}$, consider the nowhere dense real analytic subsets
$$
B_j = \left\{ (x_1, \cdots, x_{k_f^{\circ}}) : \prod_{a=1}^{k_f^{\circ}} x_a = j \right\}
$$
for nonzero $j\in \Z \cap [-2^{k_f^{\circ}}, 2^{k_f^{\circ}} ]$, $B_0=\{(0, \cdots, 0)\}$, and let $B$ be their union. 

If $f$ satisfies (RST), then the set $A(f)$ in \S \ref{sec-sst} is equidistributed in $\varphi(\mathbf{x}) d\mu_L(\mathbf{x})$,
where $\varphi:[-2, 2]^{k_f^{\circ}}\Map \R_{\ge 0}$ is a continuous function and $d\mu_L$ is the Lebesgue measure.
Therefore the set of primes $\mathfrak{p}$ of $\widetilde{F^{\circ}}$  such that
$$
\left( \frac{i_1 (a_{\mathfrak{p}})}{\sqrt{\mathbb{N}\mathfrak{p}}}, \cdots, \frac{i_{k_f^{\circ}} (a_{\mathfrak{p}})}{\sqrt{\mathbb{N}\mathfrak{p}}} \right) \in B
$$
has density $0$, since $\int_B \varphi d\mu_L = 0$.  

If $f$ satisfies ($t$-ST'), then there exists a sequence $1\le j_1 < \cdots < j_t \le k_f^{\circ}$ such that $\mathrm{pr}_{j_1, \cdots, j_t} (A(f))$
is equidistributed in the $t$-fold product of the Sato-Tate measure.  Then the set
$$
C = \left\{ (x_1, \cdots, x_{k_f^{\circ}}): \prod_{a=1}^t |x_{j_a}| < 2^{-(k_f^{\circ}-t)} \right\}
$$
is disjoint from $B_j$ for all nonzero $j\in \Z \cap [-2^{k_f^{\circ}}, 2^{k_f^{\circ}} ]$.  Since the set of $\mathfrak{p}$ such that $a_{\mathfrak{p}}=0$ has density $0$ as noted above,
the lower natural density of the set of primes $\mathfrak{p}$ such that
$
\displaystyle{\left( \frac{i_1 (a_{\mathfrak{p}})}{\sqrt{\mathbb{N}\mathfrak{p}}}, \cdots, \frac{i_{k_f^{\circ}} (a_{\mathfrak{p}})}{\sqrt{\mathbb{N}\mathfrak{p}}} \right)} \not \in B
$
is bounded from below by
\begin{equation}
\label{eqn-const}
c(k_f^{\circ}, t) := \int_{|y_1| \cdots |y_{t}| < 2^{t-k_f^{\circ}} } d\mu_{ST}(y_1) \cdots d\mu_{ST}(y_t) > 0 \,\,\,\, \mbox{ for } 1\le t<k_f^{\circ}
\end{equation}
where $d\mu_{ST}(y) = \frac{1}{2\pi} \sqrt{4-y^2}dy$ concentrated in $[-2, 2]$.

\vspace{12 pt}

\noindent (4) There is nothing to prove if $\sigma_{\widetilde{F}}(K_f)\ge 1/2$, so we assume that
$$
\lambda_0:= \lambda_{\widetilde{F}}(K_f) > k_f / 2.
$$
We return to the primes $\mathfrak{p}$ of $\widetilde{F}$ of absolute degree $1$ over $\Q$. 
The conjugacy class of $\Frob_{\mathfrak{p}}$ in $\Gal(\oQ_S / \widetilde{F})$ maps into the conjugacy class of $\Frob_p$ in $\Gal(\oQ_S/\Q)$.
(Here $\oQ_S$ is the maximal subfield of $\oQ$ unramified outside $\disc(F)\cdot\disc(K_f)\cdot\ell\cdot\mathbb{N}(\mathfrak{n})$.)
The following diagram exhibits the interaction of $\widetilde{F}$ and $K_f$ (\S \ref{sec-interplay}) and the Galois representation at hand:
$$
\xymatrix{
\Frob_{\mathfrak{p}}^{\sharp} \ar@{^{(}->}[r] \ar[d] & \Gal(\oQ_S / \widetilde{F} ) \ar[r]^-{\rho_{f,\lambda}} \ar@{^{(}->}[d] & GL_2(K_{f, \lambda}) \\
\Frob_p^{\sharp} \ar@{^{(}->}[r] & \Gal(\oQ_S / \Q) \ar[r] & \Aut (\Hom(K_f, \oQ)) \supseteq H^{\sharp}
}
$$

Let $\Gamma$ be the image of $\Gal(\oQ_S /\widetilde{F})$ in $\Aut (\Hom(K_f, \oQ))$, and let $H^{\sharp} \subseteq \Gamma$ be the nonempty subset of elements $h$ such that 
$\lambda(h, \Hom(K_f, \oQ)) = \lambda_0$; one sees easily that $H^{\sharp}$ is stable under conjugation by $\Gamma$.  
Then if $\Frob_{\mathfrak{p}}$ maps into $H^{\sharp}$,  there exists a prime ideal $\wp$ of
$K_f$ lying over $(p)$ with
$$
\dim_{\Fp} \mathcal{O}_{K_f} / \wp = \lambda_0.
$$
Now the bound in (2) prevents $\wp$ from occuring in the ideal factorisation of $a_{\mathfrak{p}}$ (\ref{eqn-factorisation}) with multiplicity $>0$.
Therefore $k(p)$ is at most the sum of the degrees of the residue class field extensions at the other primes of $K_f$ lying over $(p)$, and
$$
k(p) \le k_f - \lambda_0 = k_f \cdot \sigma_{\widetilde{F}}(K_f).
$$
The density of such primes $\mathfrak{p}$ (i.e. not dividing $\disc(F)\cdot\disc(K_f)\cdot\ell\cdot\mathbb{N}(\mathfrak{n})$, having $a_{\mathfrak{p}}\neq 0$, and whose
$\Frob_{\mathfrak{p}}$ mapping into $H^{\sharp}$) in $\widetilde{F}$ is, by the Chebotarev density theorem, equal to $|H^{\sharp}|/|\Gamma|>0$,
and we get (4).

The proof of (4bis) is parallel to that of (4), except we consider the degree-$1$ primes of $\widetilde{F^{\circ}}$ and $a_{\mathfrak{p}} \in \mathcal{O}_{K_f^{\circ}}$, and we omit it.

\vspace{12 pt}

\noindent (6) Consider a similar diagram (with $\widetilde{F}$ replaced with $\widetilde{F^{\circ}}$) to the one in (4), 
and let $\Gamma^{\circ}$ be the image of $\Gal(\oQ_S / \widetilde{F^{\circ}})$ in $\Aut(\Hom(K_f^{\circ}, \oQ))$.
This time, choose $H^{\sharp}$ to consist of those elements of $\Gamma^{\circ}$ that bisect $\Hom(K_f^{\circ}, \oQ)$. 
Again, $H^{\sharp}$ is stable under conjugation by $\Gamma^{\circ}$.

Then for any prime $\mathfrak{p}$ of $\widetilde{F^{\circ}}$ such that $\Frob_{\mathfrak{p}}$ maps into $H^{\sharp}$ and $p$ is unramified in $K_f^{\circ}$,
there are exactly $2$ primes of $K_f^{\circ}$ lying over $(p)$ with the same degree of residue class extension (namely $= k_f^{\circ}/2$).  The density of such primes $\mathfrak{p}$ in $\widetilde{F^{\circ}}$ is $|H^{\sharp}|/|\Gamma^{\circ}|>0$ by the Chebotarev density theorem.

Now the bound in (5), which is in effect because we assume (RST), keeps either of the $2$ primes from appearing in the ideal decomposition of $a_{\mathfrak{p}} \mathcal{O}_{K_f^{\circ}}$ with multiplicity $>0$,
except for a set of primes $\mathfrak{p}$ with density $0$.
\end{proof}

\begin{subremark}
The constant $c(k, t)$ for $t=1$ (where $k=k^{\circ}_f$) can be expressed:
$$
c(k, 1) = \frac{2}{\pi} \left( \frac{1}{2^{k}} \sqrt{1-\frac{1}{2^{2k}}} + \arcsin ( \frac{1}{2^k} ) \right) \,\,\,\, \mbox{ for } k\ge 2
$$
and is asymptotically $1/(\pi 2^{k-2})$ as $k\rightarrow \infty$. Here are approximate values of $c(k,t)$ for $1\le t\le k\le 6$:
\begin{center}
\begin{tabular}{c|cccccc}
 $k$ $|$ $t$ & 1 & 2 & 3 & 4 & 5 & 6 \\
 \hline
 1 & 1 & & & & & \\
 2 & 0.315 & 1 & & & & \\
 3 & 0.159 & 0.501 & 1 & & & \\
 4 & 0.0795 & 0.320 & 0.62 & 1 & & \\
 5 & 0.0398 & 0.195 & 0.45 & 0.71 & 1 & \phantom{0.0} \\
 6 & 0.0199 & 0.115 & 0.31 & 0.56 & 0.8 & 1 
\end{tabular}

\end{center}

\end{subremark}

\subsection{CM case}

\begin{subtheorem}
\label{thm-cm}
Let $f$ be a new normalised Hilbert eigencuspform of level $\mathfrak{n}\subseteq \mathcal{O}_F$ and parallel weight $(2, \cdots, 2)$.  Suppose that $f$ is of CM type (\S \ref{sec-notation}).

Denote by $M(f)$ the Andr\'e motive, whose realisations give the part of the intersection cohomology of the Hilbert modular variety
corresponding to $\{\sigma(f)\}$, where $\sigma$ ranges over all the embeddings of $K_f$ into $\oQ$. Then:

\begin{enumerate}
 \item[(1)] For all rational primes $p$ coprime to $\disc(F)\cdot \mathfrak{n}\cdot \ell$, we have
 $$
 \Newton ( \Frob_p |_{M(f)} ) \ge \HodgeP (M(f)).
 $$
 \item[(2)] For a principally abundant set of primes $p$, we have
 $$
 \Newton ( \Frob_p |_{M(f)} ) = \HodgeP (M(f)).
 $$ 
\end{enumerate}
\end{subtheorem}

\begin{proof}
Let $\lambda$ be a prime of $K_f$ lying over $(\ell)$ such that the connected component $G^{\circ}$ of the Zariski closure $G=G_{f, \lambda}$ of the image of $\rho_{f, \lambda}$ is a torus.
The argument employed in proving part (1) of Theorem \ref{thm-main} goes through without change: The non CM condition was not used.  This way we get the inequality (1)
for primes $p$ coprime to $\disc(F)\cdot \mathfrak{n}\cdot \ell$,  and for those splitting completely in $F$ and unramified in $K_f$ in addition, an integer $k(p)\in [0, k_f]$ such that
$$
 \Newton ( \Frob_p |_{M(f)}) = P(d; k_f, k(p)).
$$

Let $F^{\circ}$ be the Galois extension of $F$ cut out by the two representations with finite image:
$$
\Gal(\oQ / F^{\circ} ) = \ker ( \Gal(\oQ/F) \Map G(K_{f, \lambda}) / G^{\circ} (K_{f, \lambda}) ) \cap \ker ( \det(\rho_{f,\lambda})(1) )
$$
and let $\widetilde{F^{\circ}}$ be the compositum of $F^{\circ}$ and $\widetilde{F}$.

Since the restriction of $\rho_{f, \lambda}$ to $\widetilde{F^{\circ}}$ is then abelian and $K_f$-rational, by a theorem of Serre \cite[Chpt. III, \S 3]{serre-mcgill},
augmented with a transcendence result of Waldschmidt \cite{waldschmidt} (see Henniart \cite{henniart}), this restriction is locally algebraic (and semisimple by assumption).
Then by a theorem of Ribet \cite[\S 1.6]{ribet} (which extends that of Serre \cite[\S III.2.3]{serre-mcgill}), there exist (i) a $2$-dimensional 
$K_f$-rational vector subspace $V_0$ of $K_{f, \lambda}^{\oplus 2}$, (ii) a modulus $\mathfrak{m}$ of $\widetilde{F^{\circ}}$, and (iii) a rational representation
$$
\phi_0: S_{\mathfrak{m}} \otimes_{\Q} K_f \Map GL_{V_0}
$$
such that $\rho|_{\widetilde{F^{\circ}}}$ is isomorphic to the $\lambda$-adic representation associated with $\phi_0$.


The image of $\phi_0$ is a maximal algebraic torus of $GL_{V_0}$, since $\det\phi_0$ gives the Tate structure $\Ql(-1)$ on $\Gal(\oQ/F^{\circ})$, 
and the cyclotomic character is not divisible by $2$ as the character of any number field.

Let $K'$ be the splitting field over $K_f$ of this algebraic torus, so that $[K':K_{f}]\le 2$ and $[G_{f, \lambda}: G_{f, \lambda}^{\circ}]\le 2$.
(It is worth clarifying that unlike the $K'$ introduced in the proof of Theorem \ref{thm-main} in a similar context, this $K'$ depends only on $\rho_{f, \lambda}$ and is independent of $\mathfrak{p}$.)

For every prime $\mathfrak{p}$ of $\widetilde{F^{\circ}}$ of absolute degree $1$ and coprime to $\disc(F)\cdot\mathfrak{n}\cdot\ell$,
let $\{\alpha_{\mathfrak{p}}, \beta_{\mathfrak{p}}\} \subset K'$ be the $2$ eigenvalues of $\rho_{f, \lambda}(\Frob_{\mathfrak{p}})$.  
Since they are Weil $p$-integers by Lemma \ref{lem-Weil-number} and $\det(\rho(\Frob_{\mathfrak{p}}))=p$, we have
\begin{equation}
\label{eqn-alpha-beta-cm}
\beta_{\mathfrak{p}} = \frac{p}{\alpha_{\mathfrak{p}}} = \overline{\alpha_{\mathfrak{p}}}
\end{equation}
where the bar denotes the complex conjugation on $\Q(\alpha_{\mathfrak{p}}) = \Q(\beta_{\mathfrak{p}})$.

Now consider only those $\mathfrak{p}$ such that, in addition, $(p):= \mathfrak{p}\cap \Z$ splits completely in $K'$, a fortiori also in $\Q(\alpha_{\mathfrak{p}})$;
the resulting set is clearly principally abundant.  Then, since $(p)$ is unramified in $\Q(\alpha_{\mathfrak{p}})$, $\alpha_{\mathfrak{p}}$ cannot be totally real, and generates a CM field.
Let $\{ \wp_1, \cdots, \wp_m, \overline{\wp_1}, \cdots, \overline{\wp_m} \}$ be the set of primes of $\Q(\alpha_{\mathfrak{p}})$ lying over $(p)$, where $2m=[\Q(\alpha_{\mathfrak{p}}):\Q]$.
The equation (\ref{eqn-alpha-beta-cm}) further shows that, perhaps after renaming the primes, we get
$$
\alpha_{\mathfrak{p}}\cdot \mathcal{O}_{\Q(\alpha_{\mathfrak{p}})} = \wp_1 \cdots \wp_m \,\, \mbox{ and } 
  \,\, \beta_{\mathfrak{p}}\cdot \mathcal{O}_{\Q(\alpha_{\mathfrak{p}})} = \overline{\wp_1} \cdots \overline{\wp_m}.
$$
It follows that
$$
\Tr (\rho_{f, \lambda} (\Frob_{\mathfrak{p}})) = \alpha_{\mathfrak{p}} + \beta_{\mathfrak{p}}
$$
does not belong to any prime ideal of $\mathcal{O}_{K'}$ lying over $(p)$.  Since $\Tr(\rho_{f, \lambda} (\Frob_{\mathfrak{p}})) \in \mathcal{O}_{K_f}$, it belongs to no prime of $\mathcal{O}_{K_f}$
lying over $(p)$, either.  This proves that $k(p)=0$ and completes the proof of Theorem \ref{thm-cm}.
\end{proof}

\section{Examples}
\label{sec-numerical}

For the dimensions of and the Hecke orbits in the spaces of newforms, we rely on the information published in ``The $L$-functions and modular forms database'' \texttt{http://www.lmfdb.org/}.
We compute the slope $\sigma$ by using the polynomials given in LMFDB generating $K_f$;  sometimes the Galois group of $\widetilde{K_f} / \Q$ and the discriminant of $K_f$ are also given in LMFDB,
in which case we utilise the information also.

All the computations are for $\Gamma_0(\mathfrak{n})$ (trivial Nebentypus).  
Recall that we say two normalised eigencuspforms $f$ and $g$ with complex coefficients are conjugate (and that they belong to the same conjugacy class) if there is $\sigma\in \Aut(\C)$ such that $f^{\sigma} = g$.

\subsection{$F=\Q$}

The number of new normalised eigencuspforms $f$ of weight $2$ and level $N\le 300$ is:
$$
\sum_{N=1}^{300} \dim_{\C} S_2^{\mathrm{new}} ( \Gamma_0(N), \C ) = 2074.
$$
The degree of the field $K_f$ in this range takes the following values:
$$
k_f = [K_f: \Q] \in \{ 1, 2, 3, 4, 5, 6, 7, 8, 9, 10, 11, 12, 14, 16, 17 \}
$$
For $2070$ of the $2074$ forms $f$, parts (3) and (4) of Theorem \ref{thm-main} and Theorem \ref{thm-cm} show that $M(f)$ has an abundant set of ordinary primes. The $4$ exceptions:
\begin{itemize}
 \item There is $1$ conjugacy class of $4$ forms of level $275$ without CM, under the name 275.2.1.h in LMFDB, such that
$$
K_f = \Q (\sqrt{3}, \sqrt{11}),
$$
which is Galois with the Klein $4$-group.  There is a bisecting element, and $\sigma_{\Q}(K_f)=1/2$.

If $K_f^{\circ} \neq K_f$, then part (2) of Theorem \ref{thm-main} provides a principally abundant set of ordinary primes.  
In case $K_f^{\circ}=K_f$, part (5) and the univariate (i.e. $t=1$) Sato-Tate equidistribution (proven in \cite{taylor-sato-tate-1} and \cite{taylor-sato-tate-2})
gives an abundant set of primes $p$ (of lower density $\ge 0.0794$) such that $k(p)\in \{ 0, 1\}$; if in addition $f$ satisfies (RST), then part (6) will imply the abundance of ordinary ($k(p)=0$) primes.
\end{itemize}

\subsection{$F=\Q(\sqrt{2})$}
This quadratic field has discriminant $8$ and class number $1$.

LMFDB lists $1047$ new normalised eigencuspforms $f$ of parallel weight $(2,2)$ of level $\mathfrak{n}$ with $\mathbb{N}(\mathfrak{n})\le 350$, and the degree of $K_f$ takes the following values:
$$
k_f \in \{ 1, 2, 3, 4, 5, 6, 7, 8, 9, 11, 13 \}
$$
For $1031$ of the $1047$, parts (3) and (4) of Theorem \ref{thm-main} and Theorem \ref{thm-cm} show that $M(f)$ has an abundant set of ordinary primes. The $16$ exceptions:
\begin{enumerate}
 \item[(a)] For the $8$ forms $f$ in the classes 161.2-c and 161.3-c, we have $K_f = \Q(\sqrt{3}, \sqrt{11})$, which is Galois with the Klein $4$-group.
 As $K_f$ is linearly disjoint from $\widetilde{F}=F=\Q(\sqrt{2})$, by Proposition \ref{prop-interact} (3), $\sigma_{\widetilde{F}}(K_f)=\sigma_{\Q}(K_f)=1/2$
 and there is an element of $\Gal(\oQ/F)$ bisecting $\Hom(K_f, \oQ)$.
 
 \item[(b)] For the $8$ forms $f$ in 329.2-c and 329.3-c, $[K_f:\Q]=4$, $\widetilde{K_f}/\Q$ has group $D_8$, the image of $\Gal(\oQ/F)$ is a Klein subgroup,
 and $\sigma_{F}(K_f) = 1/2$.
\end{enumerate}
In both cases:  If $K_f^{\circ} \neq K_f$, then part (3) of Theorem \ref{thm-main} gives a principally abundant set of ordinary primes.
In case $K_f^{\circ}=K_f$, part (5) gives (unconditionally) an abundant set of primes $p$ such that $k(p)\in \{ 0, 1\}$; if $f$ satisfies (RST) in addition,
then part (6) will imply the abundance of ordinary primes.

\subsection{$F=\Q( \cos(2\pi/7))$}

This is the largest totally real subfield of the cyclotomic field $\Q(e^{2\pi i/7})$.  It is Galois over $\Q$ with group $\Z/3$ and has discriminant $49$ and class number $1$.

LMFDB lists $1075$ new normalised eigencuspforms $f$ of parallel weight $(2,2,2)$ and level $\mathfrak{n}$ with $\mathbb{N}(\mathfrak{n})\le 800$, and the degree of $K_f$ takes the following values:
$$
k_f \in \{ 1, 2, 3, 4, 5, 6, 7, 8 \}
$$
For $1048$ of the $1075$, parts (3) and (4) of Theorem \ref{thm-main} and Theorem \ref{thm-cm} show that $M(f)$ has an abundant set of ordinary primes. The $27$ exceptions:
\begin{itemize}
 \item For $f$ in the classes
 $$
 \mbox{448.1-a, 547-1-c, 547.2-c, 547.3-c, 729.1-c, 729.1-d, 743.1-a, 743.2-a and 743.3-a}
 $$
 we have $F=K_f$ and $\sigma_F(K_f)=2/3$.
 
 For each of these $f$: If $K_f^{\circ}\neq K_f$, in which case $K_f^{\circ}=\Q$, 
 then part (3) of Theorem \ref{thm-main} would give a principally abundant set of ordinary primes.
 If $K_f^{\circ}=K_f$, then the theorem only provides a principally abundant set of primes $p$ such that $k(p) \in \{ 0, 1\}$.
\end{itemize}

\subsection{$F=\Q( \cos(\pi/8))$}

This largest totally real subfield of $\Q(e^{2\pi i/16})$ has discriminant $2048=2^{11}$ and class number $1$, and is Galois over $\Q$ with group $\Z/4$. 
The nontrivial proper subgroup of $\Z/4$ allows a richer array of examples in which Theorems \ref{thm-main} and \ref{thm-cm} fall short.

LMFDB lists $6185$ new normalised eigencuspforms $f$ of parallel weight $(2,2,2,2)$ and level $\mathfrak{n}$ with $\mathbb{N}(\mathfrak{n})\le 607$, and the degree of $K_f$ takes the following values:
$$
k_f \in \{1, 2, \cdots, 12 \} \cup \{ 14, 15, 16, 17, 18, 19, 20, 22, 24, 25, 26, 27, 28, 30, 33, 39, 42 \}.
$$

For $6037$ of the $6185$, parts (3) and (4) of Theorem \ref{thm-main} and Theorem \ref{thm-cm} show that $M(f)$ has an abundant set of ordinary primes.
The $136$ confirmed exceptions and $12$ possible exceptions:
\begin{itemize}
 \item[(a)] For the $24$ forms $f$ in the classes 392.1-f, 392.2-f, 544.1-$\ell$, 544.2-$\ell$, 544.3-$\ell$ and 544.4-$\ell$, $K_f$ is Galois over $\Q$ with the Klein $4$-group and linearly disjoint from $F$ over $\Q$.
 There is a bisecting element and $\sigma_{F}(K_f)=\sigma_{\Q}(K_f)=1/2$. (cf. exceptions in $F=\Q$ and (a) in $F=\Q(\sqrt{2})$.)
 \item[(b)] For the $20$ forms $f$ in the classes 81.1-c, 289.1-f, 289.4-f, 578.1-h and 578.4-h, $k_f=4$ and $\widetilde{K_f}$ is Galois over $\Q$ with group $D_8$, so $\sigma_{\Q}(K_f)=0$.
 However, the image of $\Gal(\oQ/F)$ is a Klein $4$-group, $\sigma_{F}(K_f)=1/2$, and there is a bisecting element.  (cf. exceptions 
 (b) in $F=\Q(\sqrt{2})$.)
 \item[(c)] For the $24$ forms $f$ in 289.7-k, 289.8-k, 289.9-k and 289.10-k, $K_f$ is Galois with group $\Z/6\Z$, and $\sigma_{\Q}(K_f)=0$.
 However, the image of $\Gal(\oQ/F)$ is the subgroup $2\Z/6\Z$, $\sigma_{F}(K_f) = 1/2$, and there is a bisecting element. 
 \item[(d)] For the $68$ forms $f$ in $17$ classes (8 in level norm 289, 1 in 324 and 8 in 578), we have $K_f = F$.  Thus $\sigma_F(K_f)=3/4$ and there is \emph{no} bisecting element.
 (cf. exceptions in $F=\Q(\cos(2\pi/7))$.)
 \item[(e)] For the $12$ forms $f$ in 392.1-g and 392.2-g, $K_f$ has degree $6$ but is not cyclic over $\Q$ (hence qualitatively different from (c)).
 So far we have observed: $\sigma_{\Q}(K_f)=0$, $\sigma_F(K_f) \le 1/2$, and there is a bisecting element.
\end{itemize}
For the $f$ in (d), part (5) of Theorem \ref{thm-main} provides an abundant set of primes $p$ such that $k(p) \in \{ 0, 1\}$ unconditionally.

In the remaining cases, we have $k(p) \le k_f/2$ (resp. $k(p)<k_f/2$) for a principally abundant (resp. abundant) set of primes $p$ by part (2) (resp. by part (5)), unconditionally.
If $f$ satisfies (RST) in addition, then part (6) will provide an abundant set of ordinary ($k(p)=0$) primes.

\section{General motivic coefficients}
\label{sec-gen-mot-coeffs}

\subsection{Conjectures in a general setting}

Let $X$ be a projective variety of dimension $d$ over a number field $F$, $j:U\hookrightarrow X$ the inclusion of a smooth dense open subset, and $\pi:\mathcal{Y}\Map U$ a projective smooth scheme.
For each integer $i$ and every prime number $\ell$, form the local system on $U$,
$$
\mathcal{L}^i_{\ell} = R^i \pi_{\ast} (\Ql)
$$
and the intermediate extension 
$$
\overline{\mathcal{L}}^i_{\ell} = j_{!\ast} ( \mathcal{L}^i_{\ell}[d] ) [-d].
$$

\begin{subconjecture}
\label{conj-general}
 Let the notation be as above, and let $k$ be any integer.
\begin{enumerate}
 \item[(a)] There exists a pure Grothendieck homological motive $\mathfrak{M}=\mathfrak{M}^{k,i}$ whose $\ell$-adic \'etale realisation $\mathfrak{M}_{\ell}$ is isomorphic to $H^{k} ( X \otimes_F F^s, \overline{\mathcal{L}}^i_{\ell})$ for every $\ell$.
 \footnote{In other words, there exists a projective smooth variety $Z_k$ over $F$ and an idempotent algebraic cycle (modulo homological equivalence) $\epsilon_k$ on $Z_k \times_F Z_k$
 such that $\epsilon_k \mathfrak{h}(Z_k)$ has $\ell$-adic realisation isomorphic to $H^{k}(X \otimes_F F^s, \overline{\mathcal{L}}^i_{\ell})$.}
 \item[(a')] There exists an Andr\'e motive $M=M^{k,i}$ such that $M_{\ell} \simeq H^k (X \otimes_F F^s, \overline{\mathcal{L}}^i_{\ell})$ for every $\ell$.
\end{enumerate}
 For the following statements, we assume that (a') is true.\footnote{It appears that (a') follows from the construction of pure Nori motives realising the $\ell$-adic intersection cohomology groups, due to Ivorra and Morel \cite[\S 6]{ivorra-morel}.} 
 
 Let $e$ be an idempotent endomorphism of $M$ in the category of Andr\'e motives (with $\Q$-coefficients) and let $R$ be the direct summand of $M$ cut out by $e$,
 with the $\ell$-adic \'etale realisation $R_{\ell}$.
\begin{enumerate}
 \item[(b)] The $\ell$-adic Galois representations $R_{\ell}$ form a strictly compatible system.
 
 \item[(c)] There exists a finite set $S=S(\pi, i, k, e)$ of primes of $F$ such that, for every prime $\ell$ and $\mathfrak{p}$ outside $S$ and not dividing $\ell$, we have
 $$
 \Newton (\Frob_{\mathfrak{p}}, R_{\ell} ) \ge \HodgeP ( R_{\ell} ).
 $$
 \item[(d)] For infinitely many primes $\mathfrak{p}$ of $F$ and every prime number $\ell$, we have
 $$
 \Newton (\Frob_{\mathfrak{p}}, R_{\ell} ) = \HodgeP (  R_{\ell} ).
 $$
 \end{enumerate}
\end{subconjecture}

We note that, by Corollary \ref{cor-indep-andre}, the Hodge-Tate polygon of $R_{\ell}$ at $\lambda$ on the right hand sides is independent of the $\ell$-adic place $\lambda$ of $F$.

\begin{subproposition}
\label{prop-katz-gen-coeff}
 Assume that part (a) of Conjecture \ref{conj-general} is true and let $M$ be the Andr\'e motive of $\mathfrak{M}$,
 and that the idempotent $e$ is an \emph{algebraic} cycle, and let $\mathfrak{R}$ be the Grothendieck motive cut out by $e$ from $\mathfrak{M}$.  Then:
\begin{enumerate}
 \item[(1)] Parts (b) and (c) of the conjecture are also true for the Andr\'e motive $R$ of $\mathfrak{R}$.
 \item[(2)] If, in addition, there exists a finite extension $F'$ of $F$ such that $R_{\ell}$ restricts to an \emph{abelian} Galois representation of $F'$ for some (equivalently every) prime $\ell$,
 then for a principally abundant set of primes $\mathfrak{p}$ of $F$, we have
 $$
 \Newton (\Frob_{\mathfrak{p}}, R_{\ell} ) = \HodgeP (  R_{\ell} ),
 $$
 and, in particular, part (d) of the conjecture is also true for $R$.
\end{enumerate}
\end{subproposition}
\begin{proof}
 (1) The key point is that under the assumptions, we can use the crystalline realisation to compute the two polygons in part (c).
 Namely, for almost all $\mathfrak{p}$, we have the free $W(k(\mathfrak{p}))$-module $\mathfrak{R}_{\cris, \mathfrak{p}}$, 
 equipped with the Hodge filtration and the crystalline Frobenius $\phi_{\cris, \mathfrak{p}}$ (induced from those on $\mathfrak{M}$).

 Then, on the one hand, by Katz and Messing \cite[Th. 2]{katz-messing}, $\phi_{\cris, \mathfrak{p}}^{[k(\mathfrak{p}): \Fp]}$ has the same (multiset of) eigenvalues 
 as the $\ell$-adic Frobenius $\Frob_{\mathfrak{p}}$ on ${\mathfrak{R}}_{\ell}$.  Therefore they have the same Newton polygons.  This also proves (b).
 
 On the other hand, by Corollary \ref{cor-indep-andre}, the Hodge-Tate polygon of $R_{\ell} = \mathfrak{R}_{\ell}$ also coincides with the Hodge polygon of $\mathfrak{R}_{\cris, \mathfrak{p}}$,
 which by definition is equal to the Hodge polygon of the de Rham realisation $\mathfrak{R}_{dR}$.
 
 Now the statement (c) follows from Mazur's theorem \cite{mazur} applied to $\mathfrak{R}_{\cris, \mathfrak{p}}$.  In summary:
$$
\Newton (\Frob_{\mathfrak{p}}, R_{\ell}) = \Newton (\phi_{\cris, \mathfrak{p}}|_{\mathfrak{R}_{\cris, \mathfrak{p}}}) 
  \ge HP(\mathfrak{R}_{\cris, \mathfrak{p}}) = HP ( \mathfrak{R}_{dR} ) = \HodgeP ( R_{\ell} ).
$$

For (2), let $\widetilde{F'}$ be the normal closure of $F'$ over $\Q$ and replace $F$ with $\widetilde{F'}$, so that $\rho_{\ell}$ is abelian.
Since $\rho_{\ell}$ is $\Q$-rational and Hodge-Tate, by a theorem of Serre \cite[\S III.2.3]{serre-mcgill}, it is associated with a $\Q$-rational representation
$\phi_0: S_{\mathfrak{m}} \Map GL_{V_0}$, where $V_0$ is a $\Q$-form of $R_{\ell}$ and $\mathfrak{m}$ is a modulus of $F$.

The restriction of $\rho_{0}$ to $T_{\mathfrak{m}} \subseteq S_{\mathfrak{m}}$ can then be diagonalised: $\rho_{0}|_{T_{\mathfrak{m}}} \otimes \oQ = \chi_1 \oplus \cdots \oplus \chi_N$ and
$$
\chi_i  = \sum_{[\sigma]} n_{\sigma}(i) [\sigma],
$$
where $[\sigma]$ ranges over the characters of $T_{\mathfrak{m}}$ arising from the embeddings $\sigma$ of $F$ into $\oQ$.

Since we already know part (b) of the conjecture,  we can choose a rational prime $\ell$ that splits completely in $F$.
Then we can identify the $\sigma$ with the embeddings of $F$ into $\Ql$, once an embedding ${\iota}_{\ell}: \oQ \Map \oQ_{\ell}$ has been fixed.
Under this identification, the multiset of the Hodge-Tate weights of $\rho_{\ell}$ at any $\ell$-adic place $\lambda:F\Map \oQ_{\ell}$ is
$\displaystyle{\left\{ n_{\sigma_{0, \lambda} (i) } \right\}_{i=1, \cdots, N}}$, where $\sigma_{0, \lambda}$ is the (unique) embedding such that $\iota_{\ell}\circ\sigma_{0, \lambda} = \lambda$,
by \cite[Prop. 2, \S III.1.1]{serre-mcgill}.  

On the other hand, let $\mathfrak{p}$ be any prime of $F$ lying over any rational prime $p\neq \ell$ that splits completely in $F$.
We can also identify the embeddings $\sigma:F\Map \oQ$ with the embeddings into $\Qp$, once we fix a $p$-adic place of $\iota_p:\oQ \Map \oQ_p$.
Then the multiset of the $p$-adic valuations of $\rho_{\ell}(\Frob_{\mathfrak{p}})$ is given by $\displaystyle{\left\{ n_{\sigma_{1, \mathfrak{p}} (i)} \right\}_{i=1, \cdots, N}}$,
where $\iota_{p}\circ\sigma_{1, \mathfrak{p}}$ is the $p$-adic place $\mathfrak{p}$: See \cite[Cor. 2., \S II.3.4]{serre-mcgill}.

Since $\rho_0$ is $\Q$-rational, the two multisets are independent of $\lambda|\ell$ and $\mathfrak{p}|p$, respectively, and are equal to each other.  This proves (2).
\end{proof}

\subsection{Hilbert modular forms of motivic weights}

Let us specialise to the Baily-Borel compactification $X$ of the Hilbert modular variety $U$ defined over $\Q$ (of some level $\mathfrak{n}$) for the totally real field $F$.
For $\mathcal{Y}$, we take the universal abelian scheme $\pi:\mathcal{A}\Map U$ and the fibred product $\mathcal{A} \times_U \cdots \times_U \mathcal{A}$ over $U$.

Recall that a motivic weight $k=(k(\tau))_{\tau: F \Map \R}$ is a collection of integers $k(\tau)\ge 2$ of the same parity, for each real embedding $\tau$ of $F$.
\begin{subproposition}
 Let $f$ be a new cusp form of any motivic weight $k \neq (2, \cdots, 2)$.  Then 
\begin{enumerate}
 \item[(1)] The part $M(f)$ of the intersection cohomology of $X$ cut out by all the conjugates of $f$ satisfies parts (a), (b) and (c) of Conjecture \ref{conj-general}.
 \item[(2)] If, in addition, $f$ is of CM type, then $M(f)$ satisfies part (d) of Conjecture \ref{conj-general} also.
\end{enumerate}
\end{subproposition}
\begin{proof}
The first part follows immediately from Proposition \ref{prop-katz-gen-coeff} and the motivic construction of Galois representations 
(see Blasius-Rogawski \cite{blasius-rogawski} and the references therein).  The second part follows from Proposition \ref{prop-katz-gen-coeff}.
\end{proof}

In case $\mathcal{Y}=\mathcal{A}$ is the universal abelian scheme, the cohomology decomposes into the parts cut out by $f$ of parallel weight $(3, \cdots, 3)$.

\begin{subdefinition}
Let $G$ be a group acting on a finite set $X$.  For $g\in G$, we define $\lambda'(g, X)$ to be the \emph{smallest} of the cardinalities of the $g$-orbits in $X$;
we denote by $\lambda'(G, X)$ the supremum of $\lambda'(g, X)$ as $g$ ranges over $G$.

Given two number fields $F$ and $K$, we define
$$
\sigma'_F(K) := 1- \frac{\lambda'(\Gal(\oQ / F), \Hom (K, \oQ))}{[K:\Q]} \in \Q \cap [0, 1]
$$
\end{subdefinition}
The proof of the following is similar to that of Proposition \ref{prop-semistability}, and is omitted:
\begin{subproposition}
If $K'$ is a subfield of $K$, then $\sigma'_F(K') \le \sigma'_F(K)$.
\end{subproposition}

\begin{subdefinition}
\label{dfn-P'dki}
Let $d\ge 1$, $k\ge 1$ and $i\in [0, k]$ be integers.  We define the multiset (and the corresponding Newton polygon):
$$
P'(d; k, i) := \left( \{ 0, 2 \} ^{\otimes d} \right)^{\oplus (k-i)} \oplus \left( \{ 1, 1 \} ^{\otimes d} \right)^{\oplus i}
$$
\end{subdefinition}
This is the polygon obtained by vertically stretching $P(d; k, i)$ by a factor of $2$.

\begin{subproposition}
Let $f$ be a new normalised Hilbert eigencuspform of level $\mathfrak{n}\subseteq \mathcal{O}_F$ and parallel weight $(3, \cdots, 3)$.  
Assume that $f$ is not of CM type, and denote by $\widetilde{F}$, $F^{\circ}$, $\widetilde{F^{\circ}}$, $K_f$ and $K_f^{\circ}$
the number fields defined in the manner of \S \ref{sec-notation} and \S \ref{sec-ribet}; $\lambda$ is a prime of $K_f$ lying over a rational prime $\ell$ splitting completely in $K_f$.

Denote by $M(f)$ the Andr\'e motive (see Remark \ref{rmk-hecke-andre-motives}), whose realisations give the part of the intersection cohomology of the Hilbert modular variety
corresponding to $\{\sigma(f)\}_{\sigma}$, where $\sigma$ ranges over all the embeddings of $K_f$ into $\oQ$.

\begin{enumerate}
 \item[(1)] For all rational primes $p$ that splits completely in $F$ (equivalently in $\widetilde{F}$) and $p$ is unramified in $K_f$, then there exists an integer $k(p)\in[0, k_f]$ such that
 $$
 \Newton ( \Frob_p |_{M(f)}) = P'(d; k_f, k(p)).
 $$
 (Here $k_f = [K_f:\Q]$ and we refer to Definition \ref{dfn-Pdki} for the right hand side.)
 
 For the following parts, we only consider the primes splitting completely in $F$ and unramified in $K_f$.
 
 \item[(2)] For a principally abundant set of primes $p$, we have $k(p) \le k_f - (k_f / k_f^{\circ})$.
 \item[(3)] For an abundant set of primes $p$, we have $k(p) \le k_f \cdot \mathrm{min} ( \sigma'_{\widetilde{F}}(K_f), \sigma'_{\widetilde{F^{\circ}}} (K_f^{\circ}))$.
\end{enumerate}
\end{subproposition}

\begin{proof} (1) The proof is similar to that of part (1) of Theorem \ref{thm-main}.  The only difference is that the linear (resp. constant) coefficient of the polynomial (cf. (\ref{eqn-quadratic}))
$$
X^2 - \mathrm{Tr}(\rho(\Frob_{\mathfrak{p}})) X + \det(\rho(\Frob_{\mathfrak{p}})),
$$
has $p$-adic valuation $=2$ (resp. an integer $\ge 0$ or $\infty$) for those $p$ considered.

(2) We find a set of primes $\mathfrak{p}$ of $\widetilde{F^{\circ}}$ (resp. of $\widetilde{F}$) of density $=1$ such that 
$a_{\mathfrak{p}} = \mathrm{Tr}(\rho_{\lambda}(\Frob_{\mathfrak{p}}))$ is not divisible by at least one $p$-adic prime $\wp$ of $\widetilde{F^{\circ}}$ (resp. of $\widetilde{F}$), where $(p)=\mathfrak{p}\cap \Z$.

If $a_{\mathfrak{p}}$ is divisible by all the $p$-adic places, then it belongs to $p\cdot \mathcal{O}_{K_f^{\circ}}$ (resp. to $p\cdot \mathcal{O}_{K_f}$.
Since $M(f)$ has pure motivic weight $2$, for any archimedean place $v|\infty$, we have $|a_{\mathfrak{p}}|_v \le 2p$, hence the algebraic integer $|a_{\mathfrak{p}}/p|_{v} \le 2$,
and we form the finite set
$$
S = \{ \alpha^2 : |\alpha|_v \le 2 \, \mbox{ for all archimedean } v \, \mbox{ of } \, K_f^{\circ} \mbox{ (resp. } K_f  \mbox{)} \}.
$$

By the assumption that $f$ is not of CM type, the connected algebraic monodromy group $G^{\circ}_{f, \lambda}$ is the full $GL_2$ over $(K_f)_{\lambda})$,
and the regular map of algebraic varieties
$$
\Tr\left( \rho_{f, \lambda^{\circ}}^{\otimes 2} \otimes \det(\rho_{f,\lambda^{\circ}})^{-1} \right) : GL_2 \Map \mathbb{A}^1
$$
(resp. with $\lambda^{\circ}$ replaced with $\lambda$) is nonconstant, as $\rho_{f, \lambda^{\circ}}$ has Hodge-Tate weights both $0$ and $2$ at all $\ell$-adic places.  
Therefore the inverse image of $S\subseteq (K^{\circ}_f)_{\lambda}$ (resp. $S\subseteq (K_f)_{\lambda}$) has Haar measure $0$.  This proves the density $=1$ statement.

(3) With (2), we can now proceed with the Chebotarev type argument, similar to the one in parts (4) and (4bis) of Theorem \ref{thm-main}.  We omit the details.
\end{proof}

\acknowledgements{
The author has benefited from discussions with many mathematicians.  He thanks particularly R. Taylor for his comments on the Sato-Tate equidistributions (known and conjectural)
and his suggestion of using the uni- and bivariate distribution (which are known or accessible) in order to go beyond the Weil-Ramanujan-Petersson (``square-root'') bound.

He also thanks J.-P. Serre for pointing out that our (SST) fits within his theoretical framework developed in \cite{serre-nxp}; M. Harris for his comments on the bivariate Sato-Tate conjecture;
S. Morel and Y. Andr\'e for discussions around Andr\'e motives; L. Illusie and F. Orgogozo about constructibility theorems; R. Boltje for a discussion about tensor inductions; and N. Katz and B. Mazur for various suggestions.

Finally, he thanks the anonymous referees for their suggestions for improving the exposition.

}

\vspace{1 in}

\noindent \textbf{Addendum in response to a message in October 2024 from N. Katz} 

The references to the literature regarding ordinary reductions of abelian varieties in this article (written mostly in the academic year 2017/18) were incomplete.  I missed the following work :

\noindent W. Sawin, \textit{Ordinary primes for Abelian surfaces,} Comptes Rendus Math\'emathique,  Volume 354, Issue 6, June 2016, pages 566--568,

\noindent which would have been listed in the paragraph marked [$\ast\ast$] on the first page.

I thank N. Katz for pointing this out, and make this necessary addendum.

\vspace{0.5 in}

\bibliographystyle{amsalpha}

\begin{thebibliography}{BLGHT11}

\bibitem[A96]{andre-motifs}
Y.~Andr{\'e}.
\newblock Pour une th{\'e}orie inconditionnelle des motifs.
\newblock {\em Inst. Hautes {\'E}tudes Sci. Publ. Math.}, (83):5--49, 1996.

\bibitem[AMRT]{AMRT}
A.~Mumford D. Rapoport~M. Ash and Y.~Tai.
\newblock {\em {Smooth compactification of locally symmetric varieties}},
  volume~IV of {\em {Lie Groups: History, Frontiers and Applications}}.
\newblock Math. Sci. Press, Brookline, Mass., first edition edition, 1975.

\bibitem[BLGG11]{sato-tate-hilbert}
T.~Barnet-Lamb, T.~Gee, and D.~Geraghty.
\newblock The {S}ato-{T}ate conjecture for {H}ilbert modular forms.
\newblock {\em J. Amer. Math. Soc.}, 24(2):411--469, 2011.

\bibitem[BLGHT11]{taylor-sato-tate-2}
T.~Barnet-Lamb, D.~Geraghty, M.~Harris, and R.~Taylor.
\newblock A family of {C}alabi-{Y}au varieties and potential automorphy {II}.
\newblock {\em Publ. Res. Inst. Math. Sci.}, 47(1):29--98, 2011.

\bibitem[BBD]{bbd}
A.~A. Be\u{\i}linson, J.~Bernstein, and P.~Deligne.
\newblock {Faisceaux pervers}.
\newblock In {\em {Analysis and topology on singular spaces, {I} ({L}uminy,
  1981)}}, volume 100 of {\em {Ast{\'e}risque}}, pages 5--171. Soc. Math.
  France, Paris, 1982.

\bibitem[BR93]{blasius-rogawski}
D.~Blasius and J.~Rogawski.
\newblock Motives for {H}ilbert modular forms.
\newblock {\em Invent. Math.}, 114(1):55--87, 1993.

\bibitem[B06]{blasius-ramanujan}
Don Blasius.
\newblock Hilbert modular forms and the {R}amanujan conjecture.
\newblock In {\em Noncommutative geometry and number theory}, Aspects Math.,
  E37, pages 35--56. Friedr. Vieweg, Wiesbaden, 2006.

\bibitem[BL84]{BL}
J.-L. Brylinski and J.-P. Labesse.
\newblock {Cohomologie d'intersection et fonctions {$L$} de certaines
  vari{\'e}t{\'e}s de {S}himura}.
\newblock {\em Ann. Sci. {\'E}cole Norm. Sup. (4)}, 17(3):361--412, 1984.

\bibitem[C86]{carayol}
H.~Carayol.
\newblock Sur les repr\'{e}sentations {$l$}-adiques associ\'{e}es aux formes
  modulaires de {H}ilbert.
\newblock {\em Ann. Sci. \'{E}cole Norm. Sup. (4)}, 19(3):409--468, 1986.

\bibitem[CR]{curtis-reiner}
C.~Curtis and I.~Reiner.
\newblock {\em Methods of representation theory. {V}ol. {II}}.
\newblock Pure and Applied Mathematics (New York). John Wiley \& Sons, Inc.,
  New York, 1987.
\newblock With applications to finite groups and orders, A Wiley-Interscience
  Publication.

\bibitem[dC12]{decataldo}
M.~A. de~Cataldo.
\newblock The perverse filtration and the {L}efschetz hyperplane theorem, {II}.
\newblock {\em J. Algebraic Geom.}, 21(2):305--345, 2012.

\bibitem[dCM15]{decataldo-migliorini}
M.A. de~Cataldo and L.~Migliorini.
\newblock The projectors of the decomposition theorem are motivated.
\newblock {\em Math. Res. Lett.}, 22(4):1061--1088, 2015.

\bibitem[D71]{deligne-bourbaki}
P.~Deligne.
\newblock Formes modulaires et repr{\'e}sentations {$\ell$}-adiques.
\newblock In {\em S{\'e}minaire {B}ourbaki. {V}ol. 1968/69: {E}xpos{\'e}s
  347--363}, volume 175 of {\em Lecture Notes in Math.}, pages Exp.\ No.\ 355,
  139--172. Springer, Berlin, 1971.

\bibitem[D74]{deligne-weil}
P.~Deligne.
\newblock La conjecture de {W}eil. {I}.
\newblock {\em Inst. Hautes {\'E}tudes Sci. Publ. Math.}, (43):273--307, 1974.

\bibitem[Di13]{dimitrov-ajm}
M.~Dimitrov.
\newblock Automorphic symbols, {$p$}-adic {$L$}-functions and ordinary
  cohomology of {H}ilbert modular varieties.
\newblock {\em Amer. J. Math.}, 135(4):1117--1155, 2013.

\bibitem[EPW06]{emerton-pollack-weston}
M.~Emerton, R.~Pollack, and T.~Weston.
\newblock Variation of {I}wasawa invariants in {H}ida families.
\newblock {\em Invent. Math.}, 163(3):523--580, 2006.

\bibitem[Fa88]{faltings}
G.~Faltings.
\newblock {$p$}-adic {H}odge theory.
\newblock {\em J. Amer. Math. Soc.}, 1(1):255--299, 1988.

\bibitem[Fa89]{faltings-jami}
G.~Faltings.
\newblock Crystalline cohomology and {$p$}-adic {G}alois-representations.
\newblock In {\em Algebraic analysis, geometry, and number theory ({B}altimore,
  {MD}, 1988)}, pages 25--80. Johns Hopkins Univ. Press, Baltimore, MD, 1989.

\bibitem[Fu00]{gabber}
K.~Fujiwara.
\newblock Independence of {$\ell$} for intersection cohomology (after
  {G}abber).
\newblock In {\em Algebraic geometry 2000, {A}zumino ({H}otaka)}, volume~36 of
  {\em Adv. Stud. Pure Math.}, pages 145--151. Math. Soc. Japan, Tokyo, 2002.

\bibitem[Ha09]{harris}
M.~Harris.
\newblock Potential automorphy of odd-dimensional symmetric powers of elliptic
  curves and applications.
\newblock In {\em Algebra, arithmetic, and geometry: in honor of {Y}u. {I}.
  {M}anin. {V}ol. {II}}, volume 270 of {\em Progr. Math.}, pages 1--21.
  Birkh{\"a}user Boston, Inc., Boston, MA, 2009.

\bibitem[HSBT10]{taylor-sato-tate-1}
M.~Harris, N.~Shepherd-Barron, and R.~Taylor.
\newblock A family of {C}alabi-{Y}au varieties and potential automorphy.
\newblock {\em Ann. of Math. (2)}, 171(2):779--813, 2010.

\bibitem[He82]{henniart}
G.~Henniart.
\newblock Repr{\'e}sentations {$\ell$}-adiques ab{\'e}liennes.
\newblock In {\em Seminar on {N}umber {T}heory, {P}aris 1980-81 ({P}aris,
  1980/1981)}, volume~22 of {\em Progr. Math.}, pages 107--126. Birkh{\"a}user
  Boston, Boston, MA, 1982.

\bibitem[IM19]{ivorra-morel}
F.~Ivorra and S.~Morel.
\newblock The four operations on perverse motives (arxiv:1901.02096v1).

\bibitem[K71]{katz-ax}
N.M. Katz.
\newblock On a theorem of {A}x.
\newblock {\em Amer. J. Math.}, 93:485--499, 1971.

\bibitem[KL85]{katz-laumon}
N.M. Katz and G.~Laumon.
\newblock Transformation de {F}ourier et majoration de sommes exponentielles.
\newblock {\em Inst. Hautes {\'E}tudes Sci. Publ. Math.}, (62):361--418, 1985.

\bibitem[KM74]{katz-messing}
N.M. Katz and W.~Messing.
\newblock Some consequences of the {R}iemann hypothesis for varieties over
  finite fields.
\newblock {\em Invent. Math.}, 23:73--77, 1974.

\bibitem[Ma73]{mazur}
B.~Mazur.
\newblock Frobenius and the {H}odge filtration (estimates).
\newblock {\em Ann. of Math. (2)}, 98:58--95, 1973.

\bibitem[Mo08]{morel-jams}
S.~Morel.
\newblock Complexes pond\'{e}r\'{e}s sur les compactifications de
  {B}aily-{B}orel: le cas des vari\'{e}t\'{e}s de {S}iegel.
\newblock {\em J. Amer. Math. Soc.}, 21(1):23--61, 2008.

\bibitem[N01]{nekovar-parity-2}
J.~Nekov\'{a}\v{r}.
\newblock On the parity of ranks of {S}elmer groups. {II}.
\newblock {\em C. R. Acad. Sci. Paris S\'{e}r. I Math.}, 332(2):99--104, 2001.

\bibitem[Oc06]{ochiai}
T.~Ochiai.
\newblock On the two-variable {I}wasawa main conjecture.
\newblock {\em Compos. Math.}, 142(5):1157--1200, 2006.

\bibitem[Og82]{ogus-hodge}
A.~Ogus.
\newblock Hodge cycles and crystalline cohomology.
\newblock In {\em Hodge cycles, Motives and Shimura varieties ({L}{N}{M} 900)},
  pages 357--414. Springer-Verlag, 1982.

\bibitem[Oh83]{ohta}
M.~Ohta.
\newblock On the zeta function of an abelian scheme over the {S}himura curve.
\newblock {\em Japan. J. Math. (N.S.)}, 9(1):1--25, 1983.

\bibitem[Pa16]{patrikis}
S.~Patrikis.
\newblock Generalized {K}uga-{S}atake theory and rigid local systems, {II}:
  rigid {H}ecke eigensheaves.
\newblock {\em Algebra Number Theory}, 10(7):1477--1526, 2016.

\bibitem[Pi98]{pink-mtc}
R.~Pink.
\newblock {$\ell$}-adic algebraic monodromy groups, cocharacters, and the
  {M}umford-{T}ate conjecture.
\newblock {\em J. Reine Angew. Math.}, 495:187--237, 1998.

\bibitem[Ra78]{rapoport}
M.~Rapoport.
\newblock Compactifications de l'espace de modules de {H}ilbert-{B}lumenthal.
\newblock {\em Compositio Math.}, 36(3):255--335, 1978.

\bibitem[Ri76]{ribet}
K.~Ribet.
\newblock Galois action on division points of {A}belian varieties with real
  multiplications.
\newblock {\em Amer. J. Math.}, 98(3):751--804, 1976.

\bibitem[Se98]{serre-mcgill}
J.-P. Serre.
\newblock {\em Abelian {$\ell$}-adic representations and elliptic curves},
  volume~7 of {\em Research Notes in Mathematics}.
\newblock A K Peters, Ltd., Wellesley, MA, 1998.
\newblock With the collaboration of Willem Kuyk and John Labute, Revised
  reprint of the 1968 original.

\bibitem[Se12]{serre-nxp}
J.-P. Serre.
\newblock {\em Lectures on {$N_X (p)$}}, volume~11 of {\em Chapman \& Hall/CRC
  Research Notes in Mathematics}.
\newblock CRC Press, Boca Raton, FL, 2012.

\bibitem[Se13]{serre-oeuvres4}
J.-P. Serre.
\newblock {\em Oeuvres/{C}ollected papers. {IV}. 1985--1998}.
\newblock Springer Collected Works in Mathematics. Springer, Heidelberg, 2013.
\newblock Reprint of the 2000 edition [MR1730973].

\bibitem[SkUr14]{skinner-urban}
Ch. Skinner and E.~Urban.
\newblock The {I}wasawa main conjectures for {$\rm GL_2$}.
\newblock {\em Invent. Math.}, 195(1):1--277, 2014.

\bibitem[T89]{taylor-hilbert}
R.~Taylor.
\newblock On {G}alois representations associated to {H}ilbert modular forms.
\newblock {\em Invent. Math.}, 98(2):265--280, 1989.

\bibitem[T95]{taylor-hilbert-2}
R.~Taylor.
\newblock On {G}alois representations associated to {H}ilbert modular forms.
  {II}.
\newblock In {\em Elliptic curves, modular forms, \& {F}ermat's last theorem
  ({H}ong {K}ong, 1993)}, Ser. Number Theory, I, pages 185--191. Int. Press,
  Cambridge, MA, 1995.

\bibitem[T95]{taylor-icm1994}
R.~Taylor.
\newblock Representations of {G}alois groups associated to modular forms.
\newblock In {\em Proceedings of the {I}nternational {C}ongress of
  {M}athematicians, {V}ol.\ 1, 2 ({Z}{\"u}rich, 1994)}, pages 435--442.
  Birkh{\"a}user, Basel, 1995.

\bibitem[Ws81]{waldschmidt}
M.~Waldschmidt.
\newblock Transcendance et exponentielles en plusieurs variables.
\newblock {\em Invent. Math.}, 63(1):97--127, 1981.

\bibitem[Wn15]{wan}
X.~Wan.
\newblock The {I}wasawa main conjecture for {H}ilbert modular forms.
\newblock {\em Forum Math. Sigma}, 3:e18, 95, 2015.

\bibitem[Wi88]{wiles}
A.~Wiles.
\newblock On ordinary {$\lambda$}-adic representations associated to modular
  forms.
\newblock {\em Invent. Math.}, 94(3):529--573, 1988.

\end{thebibliography}

\end{document}